\pdfoutput=1

\documentclass[oneside,reqno,11pt]{amsart}

\usepackage{preamble}

\myexternaldocument[I:]{corank1_ASW_I}
\myexternaldocument[III:]{corank1_ASW_III}
\myexternaldocument[IV:]{corank1_ASW_IV}

\title{Co-rank \texorpdfstring{$1$}{1} Arithmetic Siegel--Weil II: \\ Local Archimedean}
\author{Ryan C. Chen}
\date{May 2, 2024}
\address{Department of Mathematics, Massachusetts Institute of Technology, 182 Memorial Drive, Cambridge, MA 02139, USA}
\email{rcchen@mit.edu}

\begin{document}
    
    \begin{abstract}
This is the second in a sequence of four papers, where we prove the arithmetic Siegel--Weil formula in co-rank $1$ for Kudla--Rapoport special cycles on exotic smooth integral models of unitary Shimura varieties of arbitrarily large even arithmetic dimension. 
Our arithmetic Siegel--Weil formula implies that degrees of Kudla--Rapoport arithmetic special $1$-cycles are encoded in the first derivatives of unitary Eisenstein series Fourier coefficients.

In this paper, we formulate and prove the key Archimedean local theorem. In the case of purely Archimedean intersection numbers, we also prove an Archimedean local arithmetic Siegel--Weil formula, relating Green currents of arbitrary degree and off-central derivatives of local Whittaker functions. The crucial input is a new limiting method, which is structurally parallel to our strategy at non-Archimedean places.
\end{abstract}

    \maketitle
    
    \tableofcontents
    
    \clearpage


    \section{Introduction}
        
This paper is a continuation of \cite{corank1_ASW_I.pdf} and will be followed by \cite{corank1_ASW_III.pdf,corank1_ASW_IV.pdf}. 
Up to missing some global context, the present paper may be read prior to all of our three other companion papers. We refer the reader to the introductions of \cite{corank1_ASW_I.pdf,corank1_ASW_IV.pdf} for further motivation, overview, and strategy for our four-part series of papers. The introduction in \cref{ssec:part_II:intro:arith_Siegel--Weil} is a abridged version of loc. cit., with emphasis on the Archimedean aspects.

        \subsection{Arithmetic Siegel--Weil}
        \label{ssec:part_II:intro:arith_Siegel--Weil}
            For the introduction, fix an imaginary quadratic field $F / \Q$ with ring of integers $\mc{O}_F$ and odd discriminant $\Delta$. Also fix an embedding $F \ra \C$. Given $m \in \Z_{\geq 0}$ and an even integer $n \in \Z$, we consider the (normalized) \emph{Siegel Eisenstein series}
    \begin{align}\label{equation:part_II:intro:results:Eisenstein}
    E^*(z,s)^{\circ}_n & \coloneqq \Lambda_m(s)^{\circ}_n \sum_{\begin{psmallmatrix} a & b \\ c & d \end{psmallmatrix} \in P_1(\Z) \backslash SU(m,m)(\Z)} \frac{\det(y)^{s - s_0}}{\det(c z + d)^n |\det (c z + d)|^{2(s - s_0)}}
    \end{align}
for the group
    \begin{equation}\label{equation:part_II:intro:results:U(m,m)}
    U(m,m) \coloneqq \left \{ h \in \Res_{\mc{O}_F / \Z} \GL_{2m} : {}^t \overline{h} \begin{pmatrix} 0 & 1_m \\ -1_m & 0 \end{pmatrix} h = \begin{pmatrix} 0 & 1_m \\ -1_m & 0 \end{pmatrix} \right \}
    \end{equation}
where $\Lambda_m(s)^{\circ}_n$ is the normalizing factor
    \begin{align}
    \Lambda_m(s)^{\circ}_n &\coloneqq \frac{(2 \pi)^{m(m-1)/2}}{(-2 \pi i)^{nm}} \pi^{m(- s + s_0)} |\Delta|^{m(m-1)/4 + \lfloor m / 2 \rfloor (s + s_0)} 
    \\
    & \mathrel{\phantom{\coloneqq}} \cdot \left ( \prod_{j = 0}^{m - 1} \Gamma(s - s_0 + n - j) \cdot L(2s + m - j, \eta^{j + n})  \right ). \notag
    \end{align}
In \eqref{equation:part_II:intro:results:U(m,m)}, the notation $1_m$ stands for the $m \times m$ identity matrix, we wrote $SU(m,m) \subseteq U(m,m)$ for the determinant $1$ subgroup, and we set $P_1 \coloneqq P \cap SU(m,m)$ for the Siegel parabolic $P \subseteq U(m,m)$ (consisting of $m \times m$ block upper triangular matrices). The variable $s \in \C$ is a complex parameter, we set $s_0 = (n - m)/2$, and the element $z = x + i y$ lies in Hermitian upper-half space (i.e. $x \in \mrm{Herm}_m(\R)$ and $y \in \mrm{Herm}_m(\R)_{>0}$; the latter means that $y$ is positive definite).\footnote{Here, the notation $\mrm{Herm}_m$ denotes a scheme over $\Spec \Z$, e.g. $\mrm{Herm}_m(\R)$ denotes $m \times m$ complex Hermitian matrices, and $\mrm{Herm}_m(\Q)$ denotes $m \times m$ Hermitian matrices with entries in $F$..} The symbol $\eta$ denotes the quadratic character associated to $F / \Q$ (via class field theory).
The sum in \eqref{equation:part_II:intro:results:Eisenstein} is convergent for $\mrm{Re}(s) > m / 2$ and admits meromorphic continuation to all $s \in \C$. When $m = 1$, the expression in \eqref{equation:part_II:intro:results:Eisenstein} is a classical Eisenstein series on the usual upper-half plane.

Given $T \in \mrm{Herm}_m(\Q)$, the Eisenstein series $E^*(z,s)^{\circ}_n$ has \emph{$T$-th Fourier coefficient}
    \begin{equation}\label{equation:part_II:intro:results:Fourier_coeff}
    E^*_T(y,s)^{\circ}_n \coloneqq 2^{m(m-1)/2} |\Delta|^{-m(m-1)/4} \int_{\mrm{Herm}_m(\Z) \backslash \mrm{Herm}_m(\R)} E^*(z, s)^{\circ}_n e^{-2 \pi i \mrm{tr}(Tz)} ~dx
    \end{equation}
for $z = x + i y$ in Hermitian upper-half space, where this integral is taken with respect to the Euclidean measure
on $\mrm{Herm}_m(\R)$. The integral is convergent for $\mrm{Re}(s) > m / 2$, and admits meromorphic continuation to all $s \in \C$. The preceding setup is as in \crefext{I:ssec:intro:Eisenstein}, and is taken from loc. cit. verbatim.

In our four-part sequence of papers, our overall objective is to prove the \emph{arithmetic Siegel--Weil formula}
    \begin{equation}\label{equation:part_II:intro:arith_Siegel--Weil}
    \frac{h_F}{w_F} \frac{d}{ds} \bigg|_{s = s_0} \frac{2 \Lambda_n(s-s_0)^{\circ}_n}{\kappa \Lambda_m(s)^{\circ}_n} E^*_T(y,s)^{\circ}_n \overset{?}{=} \widehat{\mrm{vol}}_{\widehat{\mc{E}}^{\vee}}([\widehat{\mc{Z}}(T)])
    \end{equation}
in co-rank $1$ when $m = n$, i.e. for $n \times n$ Hermitian matrices $T \in \mrm{Herm}_n(\Q)$ (with $F$-coefficients) of rank $n - 1$. We simultaneously prove the closely related ``near-central'' arithmetic Siegel--Weil formula (at $s = 1/2$) when $m = n - 1$ and $T \in \mrm{Herm}_m(\Q)$ is nonsingular.

In \cref{equation:part_II:intro:arith_Siegel--Weil}, we set $\kappa = 1$ (resp. $\kappa = 2$) if $m \neq n$ (resp. if $m = n$). For the right-hand side, there is a Kudla--Rapoport special cycle $\mc{Z}(T) \ra \mc{M}$ (finite unramified morphism of Deligne--Mumford stacks) associated to $T$. The notation $[\widehat{\mc{Z}}(T)]$ denotes a (in general conjectural) Kudla--Rapoport ``arithmetic special cycle class'' lying in a certain arithmetic Chow group $\arithCh^m(\mc{M})_{\Q}$. Here $\mc{M} \ra \Spec \mc{O}_F$ is a certain (Rapoport--Smithling--Zhang \cite{RSZ21}) everywhere smooth (stacky) integral model $\mc{M} \ra \Spec \mc{O}_F$ of a relative $n-1$-dimensional unitary Shimura variety associated to $G' \coloneqq \Res_{F / \Q} \G_m \times U(V)$, where $V$ is the non-degenerate $F / \Q$ Hermitian space of
signature $(n - 1, 1)$ which is split at all non-Archimedean places.  The right-hand side of \cref{equation:part_II:intro:arith_Siegel--Weil} denotes an \emph{arithmetic volume}, which is a real number ``defined'' by an arithmetic intersection product
    \begin{equation}
    \widehat{\mrm{vol}}_{\widehat{\mc{E}}^{\vee}}([\widehat{\mc{Z}}(T)]) ``\coloneqq" \widehat{\deg}([\widehat{\mc{Z}}(T)] \cdot \widehat{c}_1(\widehat{\mc{E}}^{\vee})^{n - m})
    \end{equation}
in an arithmetic Chow ring $\arithCh^*(\mc{M})_{\Q}$ (roughly in the sense of Gillet--Soul\'e \cite{GS87})
for a certain metrized tautological bundle $\widehat{\mc{E}}^{\vee}$ on $\mc{M}$.
More precisely, see \crefext{I:ssec:intro:arith_Siegel-Weil,I:ssec:intro:results,I:equation:intro:results:if_proper}.

The arithmetic Siegel--Weil formula is analogous to the classical Siegel--Weil formula, which relates special values of Eisenstein series and theta series associated to lattices. One main motivation for the arithmetic Siegel--Weil formula is the (in general, conjectural) theory of \emph{arithmetic theta lifting}, where one assembles the special cycle classes $\widehat{\mc{Z}}(T)$ into an automorphic ``arithmetic theta series''
    \begin{equation}\label{equation:part_II:intro:arith_Siegel-Weil:arith_theta}
    \widehat{\Theta} = \sum_T [\widehat{\mc{Z}}(T)] q^T
    \end{equation}
with $q^T \coloneqq e^{2 \pi i \mrm{tr}(Tz)}$ (with ``Fourier coefficients'' valued in an arithmetic Chow group) and one expects to use $\widehat{\Theta}$ to lift certain automorphic forms to classes in $\arithCh^m(\mc{M})_{\Q}$.

As formulated in \cref{equation:part_II:intro:arith_Siegel--Weil}, the arithmetic Siegel--Weil formula was possibly considered essentially previously known (at least up to a volume constant) for nonsingular $T \in \mrm{Herm}_n(\Q)$ by the local theorems \cite{Liu11,LZ22unitary,LL22II} (see discussion following \crefext{I:equation:intro:results:conjecture}). Our four-part series resolves the case where $T$ is singular of co-rank $1$. Non-Archimedean aspects of the arithmetic Siegel--Weil formula (along with combined non-Archimedean and Archimedean results) were open for corank $\geq 1$, prior to our work.

Set $\mc{M}_{\C} \coloneqq \mc{M} \times_{\Spec \mc{O}_F} \Spec \C$. For the purpose of this paper, the ``Archimedean part'' of the arithmetic intersection number in \cref{equation:part_II:intro:arith_Siegel--Weil} (with $m \leq n$ but not necessarily $m = n$) is of the form
    \begin{equation}\label{equation:part_II:intro:arith_Siegel--Weil:Arch_intersection_global}
    \int_{\mc{M}_{\C}} g_{T,y} \wedge c_1(\widehat{\mc{E}}^{\vee})^{n - m}
    \end{equation}
where $g_{T,y}$ is a certain real $(m - 1, m - 1)$-current on (the analytification of) $\mc{M}_{\C}$, and $c_1(\widehat{\mc{E}}^{\vee})$ is the Chern form of a certain dual metrized tautological bundle $\widehat{\mc{E}}^{\vee}$ on $\mc{M}_{\C}$. The current $g_{T,y}$ is allowed to vary with a parameter $y \in \mrm{Herm}_m(\R)_{>0}$, as is standard in this area. For nonsingular $T$, the ``Archimedean part'' of the derivative in \cref{equation:part_II:intro:arith_Siegel--Weil} seems to entail the derivative of an Archimedean local Whittaker function; for singular $T$, the corresponding Archimedean part might be more involved.

For the case of general co-rank, Garcia and Sankaran \cite{GS19} used the Mathai--Quillen theory of superconnections to define Green currents for special cycles, and used their currents to prove an ``Archimedean part'' of the arithmetic Siegel--Weil formula on (unitary and orthogonal) compact Shimura varieties.

However, in the setup explained above, the (stacky) complex Shimura varieties $\mc{M}_{\C} \coloneqq \mc{M} \times_{\Spec \mc{O}_F} \Spec \C$ will be non-proper; one needs to work with a CM extension of a totally real field $\neq \Q$ to obtain a compact Shimura variety as in \cite{GS19}. We hence need a new Archimedean result for our global arithmetic Siegel--Weil theorem. We will define our currents $g_{T,y}$ via star products and descent from the Hermitian symmetric domain. This is the classical approach of Kudla \cite{Kudla97a}, as considered by Liu \cite{Liu11} in the unitary case. Traditionally, star products were used to treat Green currents for nonsingular $T$ (or at least block-diagonally-nonsingular $T$). For $T \in \mrm{Herm}_n(\Q)$ singular of co-rank $1$, we will consider a new modification of the Green current \cref{ssec:part_II:Hermitian_domain:corank1_modification} and eventually prove our global arithmetic Siegel--Weil result using this modified Green current. We remark that Garcia--Sankaran needed a similar modification of their current \cite[{Definition 4.7}]{GS19}, but their modification is not \emph{linearly invariant}. By contrast, our modification is linearly invariant, i.e. we have
    \begin{equation}
    g_{T,y} = g_{{}^t \overline{\gamma} T \gamma, \gamma^{-1} g {}^t \overline{\gamma}^{-1}}
    \end{equation}
for any $\gamma \in \GL_m(\mc{O}_F)$. This is consistent with the expected automorphic behavior of $\widehat{\Theta}$ from \cref{equation:part_II:intro:arith_Siegel-Weil:arith_theta}. For more discussion on this, we refer the reader to our companion paper \crefext{III:sec:part_III:intro,III:sec:arith_cycle_classes}.

In this paper, we formulate and prove the ``Archimedean part'' of \cref{equation:part_II:intro:arith_Siegel--Weil} when $T$ is singular of rank $n - 1$. Our method of proof is completely different from that of \cite{GS19}. More importantly it hinges on a certain ``local limiting argument'' which plays strikingly similar role in our proof of the non-Archimedean aspects of corank $1$ arithmetic Siegel--Weil \crefext{I:ssec:intro:strategy_overview}. We further discuss this strategy next.

        \subsection{Archimedean local Arithmetic Siegel--Weil}
        \label{ssec:part_II:intro:Arch_local_ASW}
            The proof strategy is to ``take a limit" locally. In Figure \ref{figure:part_II:intro:strategy:local_limit_method} below, for a given place $v$ of $\Q$, we consider $T^{\flat} \in \mrm{Herm}_{n - 1}(\Q_v)$ with $\det T^{\flat} \neq 0$, and $T = \mrm{diag}(t, T^{\flat})$ for suitable nonzero $t \in \Q_v$. On the left, the limit refers to $t \ra 0$ in the $v$-adic topology (meaning the real topology if $v = \infty$).

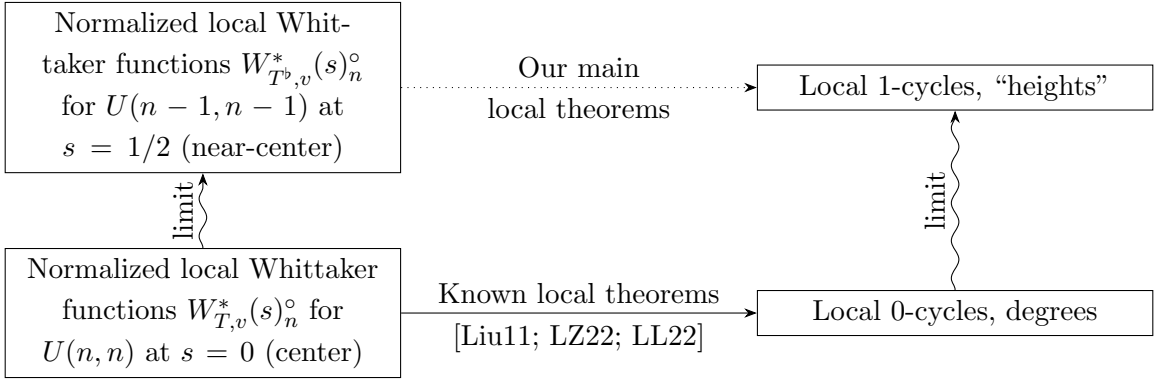
\begin{figure}[htb]
\centering
\begin{tikzpicture}[>=Stealth, node distance=3cm and 10cm, on grid, auto]
    \node (rect1) [rectangle, draw, text width=5cm, align=center] {Normalized local Whittaker functions $W^*_{T^{\flat},v}(s)^{\circ}_n$ for $U(n - 1, n - 1)$ at $s = 1/2$ (near-center)};
    \node (rect2) [rectangle, draw, text width=5cm, align=center, right=of rect1] {Local $1$-cycles, ``heights''};
    \node (rect3) [rectangle, draw, text width=5cm, align=center, below=of rect2] {Local $0$-cycles, degrees};
    \node (rect4) [rectangle, draw, text width=5cm, align=center, left=of rect3] {Normalized local Whittaker functions $W^*_{T,v}(s)^{\circ}_n$ for $U(n,n)$ at $s = 0$ (center)};

    \draw[->, dotted] (rect1) -- node[right, above] {Our main} node[right, below] {local theorems} (rect2);
    \draw[->, decorate, decoration={snake, amplitude=0.5mm, segment length=5mm}] (rect3) -- node[above, rotate=90] {limit} (rect2);
    \draw[->] (rect4) -- node[right, above] {Known local theorems} node[right, below] {\cite{Liu11,LZ22unitary,LL22II}} (rect3);
    \draw[->, decorate, decoration={snake, amplitude=0.5mm, segment length=5mm}] (rect4) -- node[above, rotate=90] {limit} (rect1);
\end{tikzpicture}
\caption{A local limiting method} \label{figure:part_II:intro:strategy:local_limit_method}
\end{figure}

In this paper, we focus on this limit argument for the Archimedean place. In this case, the lower horizontal arrow is due to Liu \cite{Liu11} (the other works \cite{LZ22unitary,LL22II} are for the non-Archimedean analogues at inert and ramified primes respectively). The preceding figure also appeared in our companion paper as \crefext{I:figure:intro:strategy:local_limit_method}. We have reproduced it here to emphasize the close analogy between our strategy at Archimedean places and our strategy at non-Archimedean places. 

In the Archimedean case, we are able to push our limit method further than at non-Archimedean places: besides $T^{\flat}$ of rank $n - 1$, we prove similar formula for $T^{\flat}$ nonsingular and non positive-definite in arbitrary rank. For the purpose of this paper, this corresponds to the arithmetic intersection number in \cref{equation:part_II:intro:arith_Siegel--Weil} being ``purely Archimedean'' (i.e. given entirely by an integral as in \cref{equation:part_II:intro:arith_Siegel--Weil:Arch_intersection_global} with no additional non-Archimedean contributions).

In the upper horizontal arrow of \cref{figure:part_II:intro:strategy:local_limit_method}, we prove an identity between a local version of \cref{equation:part_II:intro:arith_Siegel--Weil:Arch_intersection_global} and the derivative of an Archimedean local Whittaker function. For $T^{\flat} \in \mrm{Herm}_{n-1}(\Q)$ nonsingular, the relation between local Whittaker functions and the Eisenstein series Fourier coefficient is provided by the formulas
    \begin{align}
    \frac{1}{2} \frac{d}{ds}\bigg|_{s = 0} E^*_{\tilde{T}}(y,s)^{\circ}_n & = \frac{d}{d s} \bigg|_{s = 0} \left ( \frac{\Lambda_n(s)_n^{\circ}}{\Lambda_{n - 1}(s + 1/2)^{\circ}_n} E^*_{T^{\flat}}(y^{\flat}, s + 1/2)^{\circ}_n \right ) 
    \\ 
    E^*_{T^{\flat}}(y^{\flat},s)^{\circ}_n & = W^*_{T^{\flat}, \infty}(y^{\flat},s)^{\circ}_n \prod_{p} W^*_{T^{\flat}, p}(s)^{\circ}_n
    \end{align}
with $\tilde{T} = \mrm{diag}(0,T^{\flat})$ and $y = \mrm{diag}(1, y^{\flat})$; see \crefext{IV:ssec:Eisenstein:global_normalized_Fourier:singular_Fourier}. Here $W^*_{T^{\flat},\infty}(y^{\flat},s)^{\circ}_n$ (resp. $W^*_{T^{\flat}, p}(s)^{\circ}_n$) denotes a normalized Archimedean (resp. non-Archimedean) local Whittaker function; see \crefext{IV:ssec:Eisenstein:local_Whittaker:Archimedean} (resp. \crefext{IV:ssec:Eisenstein:local_Whittaker:non-Arch}).
In the Archimedean case, the local Whittaker function $W^*_{T^{\flat}}(s)^{\circ}_n$ considered in \cref{sec:part_II:local_Whittaker_Archimedean} is $W^*_{T^{\flat},\infty}(y^{\flat},s)^{\circ}_n$ evaluated at $y^{\flat} = 1$.

On the geometric side, we pass from \cref{equation:part_II:intro:arith_Siegel--Weil:Arch_intersection_global} to a local quantity via uniformization. More precisely, the complex Shimura variety $\mc{M}_{K'_f,\C}$ and its special cycles $\mc{Z}(T)_{K'_f,\C}$ admit complex uniformization (after adding sufficiently small level $K'_f \subseteq G'(\A_f)$)
    \begin{align}
    \mc{M}_{K'_f,\C}^{\mrm{an}} & \cong G'(\Q) \backslash ( \mc{D} \times G'(\A_f) / K'_f)
    \\
    \mc{Z}(T)_{K'_f} & \cong G'(\Q) \backslash \Biggl ( \coprod_{\substack{\underline{x} \in V^m \\ (\underline{x}, \underline{x}) = T}} \mc{D}(\underline{x}_{\infty}) \times \mc{D}'(\underline{x}_f) \Biggr ).
    \label{equation:part_II:intro:Arch_local_ASW:cycle_uniformization}
    \end{align}
where $\mc{D}$ is an associated Hermitian symmetric domain (parameterizing maximal negative definite subspaces in $V \otimes_{\Q} \R$), where $\mc{D}(\underline{x}_{\infty}) \subseteq \mc{D}$ is a certain ``local special cycle'' (a closed submanifold, see \cref{ssec:Hermitian_domain:local_cycles}), and $\mc{D}'(\underline{x}_f) \subseteq G'(\A_f) / K'_f$ is a subset which we call an ``away-from-$\infty$'' local special cycle. The latter is unimportant for the present paper; we omit its definition and refer to our companion paper \crefext{III:sec:Arch_uniformization} for a more detailed treatment. The disjoint union in \cref{equation:part_II:intro:Arch_local_ASW:cycle_uniformization} runs over all $m$-tuples $\underline{x}$ in $V$ with Gram matrix $T$, and the notations $\underline{x}_{\infty}$ and $\underline{x}_f$ denote the images of $\underline{x}$ in $V_\R^m \coloneqq (V \otimes_{\Q} \R)^m$ and $(V \otimes_{\Q} \A_f)^m$ respectively.

At least when $\underline{x}_{\infty}$ has nonsingular Gram matrix, the local special cycles $\mc{D}(\underline{x}_{\infty})$ have associated Green currents $[\xi(\underline{x}_{\infty})]$ on $\mc{D}$, defined by star products (\cref{ssec:Hermitian_domain:local_cycles}). Then our ``local Archimedean Siegel--Weil'' theorem is the following.

    \begin{theorem*}
    Let $\underline{x} \in V_{\R}^{m}$ be a $m$-tuple with nonsingular Gram matrix, and set $T^{\flat} = (\underline{x}, \underline{x})$. If $m \geq n - 1$ or if $T^{\flat}$ is not positive definite, we have
    \begin{equation}
    \int_{\mc{D}} [\xi(\underline{x})] \wedge c_1(\widehat{\mc{E}}^{\vee})^{n-m} = \frac{d}{d s} \bigg|_{s = -s_0} W^*_{T^{\flat}}(s)^{\circ}_n.
    \end{equation}
    where $s_0 = (n - m)/2$.
    \end{theorem*}

This appears in the text below as Theorem \ref{theorem:local_identities:Archimedean_identity:statement:main_Archimedean}. While we restricted to even $n$ for global motivational purposes (as otherwise there is some difficulty for global arithmetic Siegel--Weil on $\mc{M}$ at ramified primes), this restriction plays absolutely no role for our Archimedean local theorem (where we allow $n$ even or odd). Correspondingly, when $T^{\flat} \in \mrm{Herm}_m(\Q)$ is nonsingular and not positive-definite (giving a ``purely Archimedean'' global arithmetic intersection number), we have an analogous global arithmetic Siegel--Weil result for $\mc{M}$ of arbitrary dimension $n - 1$ \crefext{IV:theorem:part_IV:intro:results:Archimedean}. For the rest of this paper, we drop all parity assumptions on $n$.

A more precise form of our limit strategy (via \cref{figure:part_II:intro:strategy:local_limit_method}) for proving our Archimedean main local theorem will appear in \cref{ssec:Archimedean_identity:statement}. See also the sketch in \crefext{I:ssec:part_I:sketch:Archimedean}.

        \subsection{Acknowledgements}
            I thank my advisor Wei Zhang for suggesting this topic, for his dedicated support and constant enthusiasm, for insightful discussions throughout the entire course of this project, and for helpful comments on earlier drafts. I thank Tony Feng, Qiao He, Benjamin Howard, Ishan Levy, Chao Li, Keerthi Madapusi, Andreas Mihatsch, Siddarth Sankaran, Ananth Shankar, Yousheng Shi, Tonghai Yang, Shou-Wu Zhang, and Zhiyu Zhang for helpful comments or discussions.

This work was partly supported by the National Science Foundation Graduate Research Fellowship under Grant Nos. DGE-1745302 and DGE-2141064.
Parts of this work were completed at the Mathematical Sciences Research Institute (MSRI), now becoming the Simons Laufer Mathematical Sciences Institute (SLMath), and the Hausdorff Institute for Mathematics. I thank these institutes for their support and hospitality. The former is supported by the National Science Foundation (Grant No. DMS-1928930), and the latter is funded by the Deutsche Forschungsgemeinschaft (DFG, German Research Foundation) under Germany's Excellence Strategy – EXC-2047/1 – 390685813.


    \section{Hermitian symmetric domain}
        \label{sec:Hermitian_domain}
    
            \subsection{Setup}
            \label{ssec:Hermitian_domain:setup}
                We recall/fix some notation, mostly as in \cite[{\S 4B}]{Liu11} (see also \cite[{\S 2.2.2}]{GS19}). Let $n \geq 1$ be an integer, and let $V$ be the non-degenerate $\C / \R$ Hermitian space of signature $(n-1,1)$ (note that we have switched to purely local notation, in contrast with the introduction). We write $(-,-)$ for the Hermitian pairing on $V$.
Consider the Hermitian symmetric domain
    \begin{equation}
    \mc{D} = \{ \text{maximal negative definite $\C$-linear subspaces of $V$} \}.
    \end{equation}

Choosing a basis $\{e_1, \ldots, e_n\}$ of $V$ with Gram matrix $\mrm{diag}(1_{n-1}, -1)$, we take the identification
    \begin{equation}
    \begin{tikzcd}[row sep = tiny]
    \mc{D} \arrow{r}{\sim} & \{z \in \C^{n - 1} : |z| < 1 \}
    \\
    (a_1 : \cdots : a_n) \arrow[mapsto]{r} & (a_1/a_n, \ldots, a_{n-1}/a_n)
    \end{tikzcd}
    \end{equation}
and write $z_i = a_i / a_n$. Here $(a_1 : \cdots : a_n)$ stands for the complex line spanned by $a_1 e_1 + \cdots + a_n e_n$. We implicitly use the (standard) orientation $i^{n - 1} d z_1 \wedge d \overline{z}_1 \wedge \cdots \wedge d z_{n - 1} \wedge d \overline{z}_{n - 1}$ on $\mc{D}$.

We write $\mc{E}$ for the tautological line bundle over $\mc{D}$, whose fiber over a point $z \in \mc{D}$ is identified with the corresponding $\C$-line in $V$. We give $\mc{E}$ the following metric: if $w_z \in \mc{E}$ lies over $z \in \mc{D}$, set $\norm{w_z}^2 = -(w_z, w_z)$. We write $c_1(\widehat{\mc{E}})$ for the corresponding Chern form, given locally by
    \begin{equation}
    c_1(\widehat{\mc{E}}) = \frac{1}{2 \pi i} \partial \overline{\partial} \log \norm{s}^2
    \end{equation}
for local nowhere vanishing holomorphic sections $s$ of $\mc{E}$.
    
            \subsection{Local special cycles}
            \label{ssec:Hermitian_domain:local_cycles}
                Given any tuple $\underline{x} = (x_1, \ldots, x_m)$ with $x_i \in V$, there is a \emph{local special cycle}
    \begin{equation}
    \mc{D}(\underline{x}) \coloneqq \{z \in \mc{D} : z \perp x_i \text{ for all $i$} \} \subseteq \mc{D}.
    \end{equation}
This is a closed complex submanifold of $\mc{D}$.

Given $x \in V$, there is an associated global holomorphic section $s_x$ of the dual metrized tautological bundle $\widehat{\mc{E}}^{\vee}$, given by $s_x(w_z) = (x, w_z)$. For $x \in V$ and $z \in \mc{D}$, we set $R(x,z) \coloneqq \norm{s_x(z)}^2 = -(x_z, x_z)$ where $\norm{-}$ is the norm on $\widehat{\mc{E}}^{\vee}$, and $x_z$ is the orthogonal projection of $x$ to the $\C$-line $z$.

We write $\Ei(u) \coloneqq - \int_1^{\infty} e^{ut} t^{-1} ~dt$ for the exponential integral function, where $u \in \R$ is negative. We will use the asymptotics
    \begin{equation}\label{equation:Ei_function_asymptotics}
    |\Ei(u)| \leq - u^{-1} e^u \quad \quad \lim_{u \ra 0^-} (\Ei(u) - \log |u|) = \gamma,
    \end{equation}
where $\gamma$ is the Euler--Mascheroni constant. These may be verified by brief computations (omitted, but see the integral representation for $\gamma$ in \cite[{\S 12.2 Example 4}]{WW27}).

Given $x \in V$ nonzero, we set\footnote{Note that Liu instead uses $-\Ei(-2 \pi R(x,z))$ \cite[{\S 4B}]{Liu11}. This is because he considers Gram matrices $T = \frac{1}{2} (\underline{x}, \underline{x})$ while we consider Gram matrices $T = (\underline{x}, \underline{x})$ (to match our global and non-Archimedean conventions). This also affects other normalizations, e.g. our $\omega(x)$ is Liu's $\omega(\sqrt{2} x)$. \label{footnote:Green_current_normalization}.}
    \begin{equation}
    \xi(x) = -\Ei(-4 \pi R(x,z))
    \end{equation}
which is a smooth function of $z \in (\mc{D} \setminus \mc{D}(x))$ with singularity of log type along $\mc{D}(x)$ (in the sense of \cite[{(1.3.2.1)}]{GS90}). 

For locally $L^1$-forms $\xi$ on $\mc{D}$, we write $[\xi]$ for the associated current. With $x$ as above, we have the Green current equation
    \begin{equation}\label{equation:Green_current_equation_divisor}
    -\frac{1}{2 \pi i} \partial \overline{\partial} [\xi(x)] + \delta_{\mc{D}(x)} = [\omega(x)]
    \end{equation}
where $\omega(x)$ is a smooth $(1,1)$-form on $\mc{D}$ coinciding with the Kudla--Millson form up to a normalization \cite[{Proposition 4.9}]{Liu11}. Given a linearly independent tuple $\underline{x} = (x_1, \ldots, x_m) \in V^m$, we consider the current
    \begin{equation}
    [\xi(\underline{x})] \coloneqq [\xi(x_1)] * ( [\xi(x_2)] * \cdots ([\xi(x_{m-1})] * [\xi(x_m)]))
    \end{equation}
defined via star product (compare \cite[{\S 2.1.3}]{GS90}), e.g.
    \begin{equation}
    [\xi(x_1)] * ([\xi(x_2)] * [\xi(x_3)]) = \xi(x_1) \wedge \delta_{\mc{D}(x_2) \cap \mc{D}(x_3)} + \omega(x_1) \wedge \xi(x_2) \wedge \delta_{\mc{D}(x_3)} + \omega(x_1) \wedge \omega(x_2) \wedge \xi(x_3).
    \end{equation}
We then have the Green current equation
    \begin{equation}\label{equation:Green_current_equation_general}
    -\frac{1}{2 \pi i} \partial \overline{\partial} [\xi(\underline{x})] + \delta_{\mc{D}(\underline{x})} = [\omega(\underline{x})]
    \end{equation}
where $\omega(\underline{x}) \coloneqq \omega(x_1) \wedge \cdots \wedge \omega(x_m)$ (follows from \eqref{equation:Green_current_equation_general} as in the proof of \cite[{Theorem 2.4.1(i)}]{GS90}).

For any nonzero $x \in V$ and $a \in \C^{\times}$, we have
    \begin{equation}\label{equation:deform_KM_form}
    \lim_{a \ra 0} \omega(a x) = c_1(\widehat{\mc{E}}^{\vee})     
    \end{equation}
where the convergence is pointwise and uniform on compact subsets of $\mc{D} \setminus \mc{D}(x)$ (the derivatives also converge uniformly on compact subsets). This limiting statement follows upon inspecting \cite[{(2.40)}]{GS19} (see also \eqref{equation:Chern_form_taut_bundle_1} and \eqref{equation:omega(x)}). For convenience, we set $\omega(x) \coloneqq c_1(\widehat{\mc{E}}^{\vee})$ when $x = 0$.

The group $U(V)$ acts on $\mc{D}$ via the moduli description.
For any $g \in U(V)$, we have
    \begin{equation}
    g(\mc{D}(\underline{w})) = \mc{D}(g \cdot \underline{w}) \quad \quad g_* [\xi(\underline{x})] = [\xi(g \cdot \underline{x})]
    \end{equation}
where $\underline{w} \in V^m$ is any tuple and $\underline{x} \in V^m$ is any linearly independent tuple. 
    
            \subsection{Green current convergence}
            \label{ssec:Hermitian_domain:convergence}

We record some convergence estimates for the integrals appearing in our main Archimedean local identities (Section \ref{ssec:Archimedean_identity:statement}).
We work with the explicit coordinates $z = (z_1, \ldots, z_{n-1})$ on $\mc{D}$ from Section \ref{ssec:Hermitian_domain:setup} above (via the choice of basis $\{e_1, \ldots, e_n\}$ for $V$). 
For any nonzero $x \in V$, we have
    \begin{align}
    c_1(\widehat{\mc{E}}^{\vee}) = \frac{1}{2 \pi i} \partial \overline{\partial} \log R  & = \frac{1}{2 \pi i} \frac{R \partial \overline{\partial} R - \partial R \wedge \overline{\partial} R}{R^2} \label{equation:Chern_form_taut_bundle_1}
    \\
    & = \frac{1}{2 \pi i} \left ( \frac{\sum d z_j \wedge d \overline{z}_j}{1 - z \overline{z}} + \frac{(\sum \overline{z}_j d z_j) \wedge (\sum z_j d \overline{z}_j)}{(1 - z \overline{z})^2} \right ). \label{equation:Chern_form_taut_bundle_2}
    \end{align}
and
    \begin{align}\label{equation:omega(x)}
    \omega(x) = -\frac{1}{2 \pi i} \partial \overline{\partial} \xi(x) = \frac{1}{2 \pi i} e^{-4 \pi R} \left ( \frac{-4 \pi \partial R \wedge \overline{\partial} R}{R} + \frac{\partial \overline{\partial} R}{R} - \frac{\partial R \wedge \overline{\partial} R}{R^2} \right )
    \end{align}
on $\mc{D} \setminus \mc{D}(x)$, where $R \coloneqq R(x,z)$ for short.

\begin{lemma}\label{lemma:convergence:KM_form_growth_bound}
For any fixed $x \in V$ (possibly $x = 0$) with $\omega(x) = \sum_{i,j} \omega(x)_{i,j} d z_i \wedge d \overline{z}_j$, the functions $(1 - z \overline{z})^3 \omega(x)_{i,j}$ are bounded on $\mc{D}$.
\end{lemma}
\begin{proof}
If $x = \sum a_j e_j$, we have
    \begin{equation}
    R(x,z) = \frac{(a_1 \overline{z}_1 + \cdots + a_{n - 1} \overline{z}_{n - 1} - a_n)(\overline{a}_1 z_1 + \cdots + \overline{a}_{n - 1} z_{n - 1} - \overline{a}_n)}{(1 - z \overline{z})}.
    \end{equation}
This expression and the formulas for $\omega(x)$ (see above) yield the lemma via straightforward computation (omitted).
\end{proof}

\begin{lemma}\label{lemma:convergence:R_growth_boundary}
Let $\underline{x} = (x_1, \ldots, x_{m}) \in V^{m}$ be an $m$-tuple with nonsingular Gram matrix $(\underline{x}, \underline{x})$. Assume either that $m \geq n - 1$ or that $(\underline{x}, \underline{x})$ is not positive definite. Then exists $\e > 0$ such that
    \begin{equation}
    \sum_{i=1}^{n-1} R(x_i, z) > \frac{\e}{1 - z \overline{z}}
    \end{equation}
for all $z \in \mc{D}$ with $|z| \gg 0$.
\end{lemma}
\begin{proof}
Given $x  = \sum_j a_j e_j \in V$, we use the temporary notation $x \cdot z \coloneqq a_1 \overline{z}_1 + \cdots + a_{n-1} \overline{z}_{n-1} - a_n$ for $z = (z_1, \ldots, z_{n-1}) \in \C^{n-1}$. Note $R(x,z) = |x \cdot z|^2 (1 - z \overline{z})^{-1}$ for $z \in \mc{D}$. View $\C^{n - 1}$ as a standard coordinate chart in the projective space of lines in $V$ (i.e. the lines which are not orthogonal to $e_n$). The zeros of $\sum_i |x_i \cdot z|^2$ on $\C^{n-1}$ correspond to those lines in $V$ (in the given chart) which are orthogonal to $\operatorname{span}(\underline{x})$. This (closed) set of zeros is disjoint from the set $\{z \in \C^{n - 1} : |z| = 1\}$, which corresponds to isotropic lines in $V$ (i.e. no isotropic lines in $V$ are orthogonal to $\operatorname{span}(\underline{x})$). Hence $\sum_i R(x_i,z)(1 - z \overline{z})$ is bounded below (as a function of $z \in \mc{D}$) by a positive constant as $|z| \ra 1$.
\end{proof}

\begin{lemma}\label{lemma:convergence:integral_convergence}
Let $\underline{x} = (x_1, \ldots, x_{m}) \in V^{m}$ be an $m$-tuple with nonsingular Gram matrix $(\underline{x}, \underline{x})$. Assume either that $m \geq n - 1$ or that $(\underline{x}, \underline{x})$ is not positive definite.

Let $\omega = \sum \omega_{I,J} d z_I \wedge d \overline{z}_J$ (multi-indices) be any smooth complex differential form on $\mc{D}$ such that each $(1 - z \overline{z})^b \omega_{I,J}$ is bounded on $\mc{D}$ for some real constant $b \gg 0$.
Then the integral
    \begin{equation}
    \int_{\mc{D}} \xi(x_1) \omega(x_2) \wedge \cdots \wedge \omega(x_{m}) \wedge \omega
    \end{equation}
is absolutely convergent.
\end{lemma}
\begin{proof}
After making a unitary change of basis for $V$, we may assume
    \begin{equation}
    x_1 =
    \begin{cases}
    a e_n & \text{if $(x_1,x_1) < 0$} \\
    a e_1 & \text{if $(x_1,x_1) > 0$} \\
    e_{n-1} + e_n & \text{if $(x_1,x_1) = 0$}
    \end{cases}
    \end{equation}
for some nonzero $a \in \R$ (where $(e_1, \ldots, e_n)$ is the basis of $V$ used to define the coordinates $(z_1, \ldots, z_{n-1})$ in Section \ref{ssec:Hermitian_domain:setup}). This will aid calculation in coordinates.

Lemma \ref{lemma:convergence:KM_form_growth_bound} shows that it is enough to check (absolute) convergence of
    \begin{equation}\label{equation:convergence:integral_convergence:dominating_integral}
    \int_{\mc{D}} \xi(x_1) e^{-4 \pi (R(x_2, z) + \cdots + R(x_{m}, z))} (1 - z \overline{z})^{-b}
    \end{equation}
for any $b \in \R$ (for the Euclidean measure on $\mc{D}$). It is enough to check convergence when $b \gg 0$, so we assume $b \geq n$ for convenience.

Set
    \begin{equation}
    u_j \coloneqq \frac{\operatorname{Re}(z_j)}{\sqrt{1 - z \overline{z}}}
    \quad \quad
    v_j \coloneqq \frac{\operatorname{Im}(z_j)}{\sqrt{1 - z \overline{z}}}
    \end{equation}
for $j = 1, \ldots, n-1$. A change of variables gives
    \begin{align}
    & \int_{\mc{D}} \xi(x_1) e^{-4 \pi (R(x_2, z) + \cdots + R(x_{m}, z))} (1 - z \overline{z})^{-b} \label{equation:convergence_integral_new_coordinates}
    \\
    & = \int_{\R^{2(n-1)}} \xi(x_1) e^{-4 \pi (R(x_2, z) + \cdots + R(x_{m}, z))} ( 1 + |u|^2 + |v|^2)^{b - n}, \notag
    \end{align}
where $|u|^2 \coloneqq \sum_j u_j^2$ and $|v|^2 \coloneqq \sum_j v_j^2$, with $R(x_i,z)$ a function of $u,v$ via \eqref{equation:convergence_integral_new_coordinates}, and with the Euclidean measure $d u_1 ~d v_1 \cdots  d u_{n - 1} d v_{n - 1}$ understood on the right-hand side.

The asymptotics for $\Ei(u)$ as in \eqref{equation:Ei_function_asymptotics} show it is enough to check convergence of the integrals
    \begin{align}
    & &&\int_{\R^{2(n-1)}} e^{-4 \pi(R(x_1, z) + R(x_2, z) + \cdots R(x_{m}, z))} (1 + |u|^2 + |v|^2)^{b - n} \label{equation:convergence:decaying_part}
    \\
    & \text{and} && \int_{\substack{\R^{2(n-1)} \\ R(x_1, z) \leq 1/(8 \pi)}} \log(4 \pi R(x_1, z)) e^{-4 \pi (R(x_2, z) + \cdots + R(x_{m}, z))} ( 1 + |u|^2 + |v|^2)^{b - n} \label{equation:convergence:log_singular_part}
    \end{align}
(where the second integral is over the set of $(u,v) \in \R^{2(n-1)}$ satisfying $R(x_1, z) \leq 1/(8 \pi)$).

Since we have $(1 - z \overline{z})^{-1} = 1 + |u|^2 + |v|^2$, Lemma \ref{lemma:convergence:R_growth_boundary} implies that \eqref{equation:convergence:decaying_part} is absolutely convergent (by exponential decay of the integrand as $|u|^2 + |v|^2 \ra \infty$). 

For convergence of \eqref{equation:convergence:log_singular_part}, the same lemma shows that it is enough to check convergence of the integral
    \begin{equation}\label{equation:convergence_log_singular_part_2}
    \int_{\substack{\R^{2(n-1)} \\ R(x_1, z) \leq 1/(8 \pi)}} \log(8 \pi R(x_1, z)) e^{-4 \pi \e (1 + |u|^2 + |v|^2)} ( 1 + |u|^2 + |v|^2)^{b - n} 
    \end{equation}
for all $\e > 0$ (using also $R(x_1, z) \leq 1/(8 \pi)$). We check this convergence via casework.

\emph{Case when $(x_1, x_1) < 0$:} In this case, we have $R(x_1, z) = a^2 (1 + |u|^2 + |v|^2)$. The integrand in \eqref{equation:convergence:log_singular_part} is bounded on the compact set $\{ (u,v) \in \R^{2(n-1)} : R(x_1,z) \leq 1/(8 \pi) \}$, hence the integral is convergent.

\emph{Case when $(x_1, x_1) > 0$:} In this case, we have $R(x_1, z) = a^2 (u_1^2 + v_1^2)$. To check convergence of \eqref{equation:convergence_log_singular_part_2}, it is enough to check that
	\begin{equation}
	\int_{\substack{\R^{2(n-1)} \\ a^2(u_1^2 + v_1^2) \leq 1/(8 \pi)}} \log(4 \pi a^2 (u_1^2 + v_1^2)) e^{-4 \pi \e (1 + u_2^2 + v_2^2 + \cdots + u_{n-1}^2 + v_{n-1}^2)} (1 + u_2^2 + v_2^2 + \cdots + u_{n-1}^2 + v_{n-1}^2)^{b - n}
	\end{equation}
is convergent (using $R(x_1, z) \leq 1/(8 \pi)$). The integral over $(u_1,v_1)$ converges because the singularity at $u_1 = v_1 = 0$ is logarithmic, and the integral over $(u_2, v_2, \ldots, u_{n-1}, v_{n-1})$ converges because of the exponential decay.

\emph{Case when $(x_1, x_1) = 0$:} In this case, we have $R(x_1, z) = (u_1 - \sqrt{1 + |u|^2 + |v|^2})^2+ v_1^2$. Under the condition $R(x_1, z) \leq 1/(8 \pi)$, we may bound $|\log R(x_1,z)| \leq C \cdot (1 + |u_1|)$ for some constant $C > 0$. 
To check convergence of \eqref{equation:convergence_log_singular_part_2}, it is thus enough to check that
    \begin{equation}
    \int_{\R^{2(n-1)}} (1 + |u_1|) e^{-4 \pi \e (1 + |u|^2 + |v|^2)} (1 + |u|^2 + |v|^2)^{b - n}
    \end{equation}
is convergent, which follows from exponential decay of the integrand.
\end{proof}

\begin{remark}\label{remark:local_identities:Archimedean_identity:convergence:non_convergent}
The convergence result of Lemma \ref{lemma:convergence:integral_convergence} fails in general if $m < n -1$ and $(\underline{x}, \underline{x})$ is positive definite. For example, if $n = 3$, if $m = 1$, and if $x \in V$ with $(x,x) > 0$, the integral
    \begin{equation}
    \int_{\mc{D}} \xi(x) \wedge \omega(0)^2 
    \end{equation}
is not absolutely convergent.
\end{remark}

            \subsection{Co-rank \texorpdfstring{$1$}{1} modification of current}
            \label{ssec:part_II:Hermitian_domain:corank1_modification}
                
As discussed in the introduction, we propose a new \emph{linearly invariant} modification for Green currents in the co-rank $1$ singular situation. In this paper, we focus on the local version (on the Hermitian symmetric domain) but also sketch the global version (on the complex Shimura variety, constructed from the local one via uniformization). The global version (with more details and in greater generality) is treated in \crefext{III:ssec:Arch_uniformization:Archimedean}.

We switch notation from the rest of Section \ref{sec:Hermitian_domain}: let $F / \Q$ be an imaginary quadratic field and let $V$ be a signature $(n - 1, 1)$ non-degenerate $F / \Q$ Hermitian space. Fix an embedding $F \ra \C$. Form the Hermitian symmetric domain $\mc{D}$ as in Section \ref{ssec:Hermitian_domain:setup} (associated to $V_{\R}$ in the present notation). Fix a non-degenerate full-rank Hermitian $\mc{O}_F$-lattice $L \subseteq V$.

Fix $T \in \mrm{Herm}_m(\Q)$. In \cref{ssec:part_II:Hermitian_domain:corank1_modification}, we assume that $m \geq n - 1$ if $T$ is positive definite. If $T$ is singular, we assume $m = n$ and that $\rank(T) = n - 1$.
Fix $y \in \mrm{Herm}_m(\R)_{>0}$, i.e. $y$ is any positive definite complex Hermitian matrix. As is typical in Kudla's program, we allow Green 
currents to vary with the auxiliary parameter $y$.

Consider $\underline{x} \in V$ with $(\underline{x}, \underline{x}) = T$. If $T$ is nonsingular, set
    \begin{equation}\label{equation:Arch_uniformization:Archimedean_nonsingular_current}
    [\xi(\underline{x}, y)] \coloneqq [\xi(\underline{x} \cdot a)]
    \end{equation}
for a choice of $a \in \GL_m(\C)$ satisfying $a {}^t \overline{a} = y$, with $[\xi(\underline{x} \cdot a)]$ the current from Section \ref{ssec:Hermitian_domain:local_cycles}. We will not check that the current $[\xi(\underline{x},y)]$ is independent of the choice of $a$, but the intersection numbers appearing in our main results will not depend on $a$ (Remark \ref{remark:Archimedean_identity:statement:linear_invariance_geometric}, also the ``linear invariance'' from \cite[{Proposition 4.10}]{Liu11} when $m = n$).

Next, suppose that $T$ is singular, with $m = n$ and $\rank(T) = n - 1$. First consider the case when $T = \mrm{diag}(0, T^{\flat})$ where $T^{\flat}$ is nonsingular of rank $n - 1$. If $(\underline{x}, \underline{x}) = T$, we must have $\underline{x} = [0, x_2, \ldots, x_n] \in V^n$. Set $\underline{x}^{\flat} = [x_2, \ldots, x_n]$. There is a decomposition
    \begin{equation}\label{equation:Arch_uniformization:Archimedean:decomp_y_diagonal}
    y = 
    \begin{pmatrix}
    1 & c \\
    0 & 1
    \end{pmatrix}
    \begin{pmatrix}
    y^{\#} & 0 \\
    0 & y^{\flat}
    \end{pmatrix}
    \begin{pmatrix}
    1 & 0 \\
    {}^t \overline{c} & 1
    \end{pmatrix}
    \end{equation}
for uniquely determined $c \in M_{1,n-1}(\C)$, $y^{\#} \in \R_{>0}$, and $y^{\flat} \in \mrm{Herm}_{n - 1}(\R)_{>0}$. We then set
    \begin{equation}
    [\xi(\underline{x}, y)] \coloneqq c_1(\widehat{\mc{E}}^{\vee}) \wedge [\xi(\underline{x}^{\flat}, y^{\flat})] - \log(y^{\#}) \cdot \delta_{\mc{D}(\underline{x})}.
    \end{equation}
    
For $T$ not necessarily block-diagonal, we define $[\xi(\underline{x}, y)]$ by the linear invariance requirement
    \begin{equation}
    [\xi(\underline{x},y)] \coloneqq [\xi(\underline{x} \cdot \gamma^{-1}, \gamma y {}^t \overline{\g})] \mod \sum_{\substack{p \text{ such that} \\ \gamma \not \in \GL_n(\mc{O}_F \otimes_{\Z} \Z_{(p)})}} \Q \cdot \log p \cdot \delta_{\mc{D}(\underline{x})}
    \end{equation}
for all $\gamma \in \GL_n(F)$, where ${}^t \overline{\gamma}$ means conjugate transpose, and where ``$\mrm{mod}$'' means that the equality (of currents) holds up to adding an element of the displayed sum. 

Equivalently, suppose $\gamma \in \GL_n(F)$ is any element such that ${}^t \overline{\gamma}^{-1} T \gamma^{-1} = \mrm{diag}(0, T^{\flat})$ is block diagonal with $T^{\flat}$ nonsingular. Write $\underline{x} \cdot \gamma^{-1} = [0, x_1^{\flat}, \ldots, x_{n - 1}^{\flat}]$, set $\underline{x}^{\flat}_{\gamma} = [x_1^{\flat}, \ldots, x_{n - 1}^{\flat}]$, and decompose
    \begin{equation}
    \gamma y {}^t \overline{\g} = 
    \begin{pmatrix}
    1 & c \\
    0 & 1
    \end{pmatrix}
    \begin{pmatrix}
    y_{\gamma}^{\#} & 0 \\
    0 & y_{\gamma}^{\flat}
    \end{pmatrix}
    \begin{pmatrix}
    1 & 0 \\
    {}^t \overline{c} & 1
    \end{pmatrix},
    \end{equation}
as above (temporary notation). We then have
    \begin{equation}\label{equation:arch_uniformization:Archimedean:singular_diagonalized_current}
    [\xi(\underline{x},y)] = c_1(\widehat{\mc{E}}^{\vee}) \wedge [\xi(\underline{x}^{\flat}_{\gamma}, y^{\flat}_{\gamma})] - \log(\tilde{y}^{\#}) \cdot \delta_{\mc{D}(\underline{x})}
    \end{equation}
for a positive real number $\tilde{y}^{\#}$ uniquely determined by $T$ and $y$. Indeed,
we require
    \begin{equation}
    \log(\tilde{y}^{\#}) = \log(y^{\#}_{\gamma}) \quad \mrm{mod} \quad \sum_{\substack{p \text{ such that} \\ \gamma \not \in \GL_n(\mc{O}_F \otimes_{\Z} \Z_{(p)})}} \Q \cdot \log p.
    \end{equation}
For any fixed prime $p$, we can always find $\gamma \in \GL_n(\mc{O}_F \otimes_{\Z} \Z_{(p)})$ such that ${}^t \overline{\gamma}^{-1} T \gamma^{-1}$ is block diagonal as above. The preceding expression thus characterizes $\tilde{y}^{\#}$ uniquely.\footnote{For any integer $N$, set $\R_N \coloneqq \R / (\sum_{p \mid N} \Q \cdot \log p)$. For any set of integers $\{N_i\}_{i \in I}$, the diagram
    \begin{equation}
    \begin{tikzcd}[ampersand replacement=\&]
    \R_{\mrm{gcd}(\{N_i\}_{i \in I})} \arrow{r} \& \bigoplus_{i \in I} \R_{N_i} \arrow[shift left]{r} \arrow[shift right]{r} \& \bigoplus_{(i,i') \in I^2} \R_{N_i N_{i'}}
    \end{tikzcd}
    \end{equation}
is an equalizer in the category of sets. 
}
In all cases above ($T$ singular or not), note $[\xi(\underline{x},y)] = [\xi(\underline{x} \cdot \gamma^{-1}, \gamma y {}^t \overline{\gamma})]$ for all $\gamma \in \GL_m(\mc{O}_F)$ (``linear invariance'').

We quickly sketch how $[\xi(\underline{x},y)]$ descends to the (stacky) complex Shimura variety $\mc{M}_{\C} \coloneqq \mc{M} \times_{\Spec \mc{O}_F} \Spec \C$ via uniformization. A more detailed discussion (in greater generality) appears in \crefext{III:sec:Arch_uniformization}. To match the geometric setup of \crefext{I:ssec:part_I:arith_intersections:integral_models}, we take $L$ to be self-dual and $2 \nmid \Delta$ (see \crefext{III:sec:Arch_uniformization} for the general case, which is similar).

Write $\A_f$ for the finite ad\`ele ring of $\Q$. With $G' \coloneqq (\Res_{F / \Q} \G_m) \times U(V)$ and $K'_f \subseteq G'(\A_f)$ a small factorizable open compact subgroup, we consider $\mc{M}_{K'_f}$ (i.e. $\mc{M}$ with level $K'_f$ structure \crefext{III:ssec:ab_var:level_structure}) and have complex uniformization
    \begin{equation}
    G'(\Q) \backslash (\mc{D} \times G'(\A_f)/K'_f) \xra{\sim} \mc{M}_{K'_f,\C}^{\mrm{an}}
    \end{equation}
where the right-hand side denotes the analytification of $\mc{M}_{K'_f, \C}$. Assume we have a factorization $K'_f = K_{0,f} \times K_f$ (with $K_{0,f} \subseteq (\Res_{F/\Q} \G_m)(\A_f)$ and $K_f \subseteq U(V)(\A_f)$). Define the sets
    \begin{align}
    J_{\infty}(T) & \coloneqq GU(V_0)(\A_f) / K_{0,f} \times \coprod_{\substack{\underline{x} \in V^m \\ (\underline{x}, \underline{x}) = T}} \mc{D}(\underline{x}_f).
    \\
    \mc{D}(\underline{x}_f) & \coloneqq \{ g \in U(V)(\A_f) / K_f : g^{-1} x_i \in L \otimes_{\Z} \hat{\Z}^p \text{ for all $x_i \in \underline{x}_f$} \}
    \end{align}
where $\underline{x}_f \in V(\A_f)^m$ is the image of $\underline{x}^m$ under $V(\Q)^m \ra V(\A_f)^m$.
We think of $\mc{D}(\underline{x}_f)$ as an ``away-from-$\infty$'' special cycle, and it appears in complex uniformization of global special cycles \crefext{III:sec:Arch_uniformization}.

In one of our companion papers, we show that the groupoid $[G'(\Q) \backslash J_{\infty}(T)]$ has finite stabilizers and finitely many isomorphism classes (\crefext{IV:lemma:local_Siegel-Weil:uniformization_degree:main}). For each $j \in J_{\infty}(T)$, there is a corresponding map
    \begin{equation}
    \Theta_j \colon \mc{D} \ra \mc{M}^{\mrm{an}}_{K'_f, \C}
    \end{equation}
induced by the uniformization morphism $\mc{D} \times G'(\A_f) / K'_f \ra \mc{M}^{\mrm{an}}_{\C}$.

Pick any set of representatives $J \subseteq J_{\infty}(T)$ for the quotient $G'(\Q) \backslash J_{\infty}(T)$. Note that every $j \in J$ has an associated tuple $\underline{x} \in V^m$. We descend the local currents $[\xi(\underline{x},y)]$ to a current $g_{T,y}$ on $\mc{M}^{\mrm{an}}_{K'_f, \C}$ by the formula
    \begin{equation}
    g_{T,y} \coloneqq \sum_{j \in J} \frac{1}{|\Aut(j)|} \Theta_{j,*} [\xi(\underline{x},y)]
    \end{equation}
where $\Theta_{j,*}$ denotes a pushforward of currents. For singular $T$, see the convergence estimates in Section \ref{ssec:Hermitian_domain:convergence}. The current $g_{T,y}$ does not depend on the choice of $J$. Varying over levels $K'_f$ defines a current $g_{T,y}$ on $\mc{M}_{\C}$ (\crefext{III:ssec:arith_cycle_classes:horizontal}).


    \section{Archimedean local Whittaker functions}
    \label{sec:part_II:local_Whittaker_Archimedean}

        \subsection{Additional notation on \texorpdfstring{$U(m,m)$}{U(m,m)}} 
            \label{ssec:part_II:local_Whittaker_Archimedean:group}
                Recall the group $U(m,m)$ defined in \eqref{equation:part_II:intro:results:U(m,m)}. There we defined $U(m,m)$ with respect to the degree $2$ finite flat map $\Z \ra \mc{O}_F$, but we can make the same definition with respect to any finite locally free morphism $A \ra B$ of commutative rings such that $B$ is equipped with an involution $b \mapsto \overline{b}$ over $A$. In \cref{sec:part_II:local_Whittaker_Archimedean}, we are interested in the case where $A \ra B$ is $\R \ra \C$, and sometimes use the alternative notation $H \coloneqq U(m,m)$.

We consider
    \begin{align}
    P &\coloneqq  \left \{ h = \begin{pmatrix} * & * \\ 0 & * \end{pmatrix} \in H \right \}  && M \coloneqq  \left \{ m(a) = \begin{pmatrix} a & 0 \\ 0 & {}^t \overline{a}^{-1} \end{pmatrix} : a \in \Res_{B/A} \GL_m \right \} \\
    w & \coloneqq \begin{pmatrix} 0 & 1_m \\ -1_m & 0 \end{pmatrix} && N \coloneqq  \left \{ n(b) = \begin{pmatrix} 1_m & b \\ 0 & 1_m \end{pmatrix} : b \in \mrm{Herm}_m  \right \} 
    \end{align}
where $P$, $M$, and $N$ are subgroups of $H$. We have $P(R) = M(R) N(R)$ for all $A$-algebras $R$.

We consider the standard maximal compact subgroup
    \begin{equation}
    K^{\circ} \coloneqq \left \{ \begin{pmatrix} a & b \\ - b & a \end{pmatrix} \in H(\R) : a {}^t \overline{a} + b {}^t \overline{b} = 1_m \text{ and } a {}^t \overline{b} = b {}^t \overline{a} \right \} \subseteq H(\R)
    \end{equation}
We write $U(m) \subseteq \GL_m(\C)$ for (the real points of) the unitary group for the usual positive definite rank $m$ complex Hermitian space (specified by the Gram matrix $1_m$). There is an isomorphism $K^{\circ} \ra U(m) \times U(m)$ sending the displayed matrix to $(a + i b, a - i b) \in U(m) \times U(m)$ (see e.g. \cite[{\S 2.5.1}]{GS19}). We have $H(\R) = P(\R) K^{\circ}$. 
    
            \subsection{Principal series}
            \label{ssec:part_II:local_Whittaker_Archimedean:principal_series}
                Characters are assumed continuous and unitary unless specified otherwise. Form the group $H = U(m,m)$ as in Section \ref{ssec:part_II:local_Whittaker_Archimedean:group}
Given $s \in \C$ and a character $\chi \colon \C^{\times} \ra \C^{\times}$, we form the local \emph{degenerate principal series}
    \begin{equation}
    I(s,\chi) \coloneqq \Ind_{P(\R)}^{H(\R)}(\chi \cdot |-|_{F}^{s+m/2}).
    \end{equation}
This is an unnormalized induction, consisting of smooth and $K^{\circ}$-finite functions $\Phi \colon H(\R) \ra \C$ satisfying
    \begin{equation}
    \Phi(m(a) n(b) h, s) = \chi(a) |\det a|_F^{s + m/2}
    \end{equation}
for all $m(a) \in M(\R)$ and $n(b) \in N(\R)$ and $h \in H(\R)$. Here we wrote $\chi(a) \coloneqq \chi(\det a)$ for short. A section $\Phi(h,s)$ of $I(s, \chi)$ is \emph{standard} if $\Phi(k,s)$ is independent of $s$ for any fixed $k \in K^{\circ}$. 

Suppose $\chi|_{\R^{\times}} = \eta^{n}$ for some integer $n$, where $\eta \coloneqq \mrm{sgn}(-)$ is the sign character. In this case, let $k(\chi) \in \Z$ be the integer satisfying
    \begin{equation}
    \chi(z) = (z / |z|^{1/2}_{\C})^{k(\chi)},
    \end{equation}
(where $|-|_{\C}$ is the square of the usual Euclidean norm) and let $\Phi^{(n)}$ be the unique standard section of $I(s, \chi)$ of scalar weight
    \begin{equation}
    (n_1, n_2) \quad \quad \text{where} \quad \quad n_1 = \frac{n + k(\chi)}{2} \quad \text{and} \quad n_2 = \frac{-n + k(\chi)}{2}
    \end{equation}
such that $\Phi^{(n)}(1, s) = 1$ (as in \cite[{\S 3.2, \S 3.3}]{GS19}). The scalar weight condition means that
    \begin{equation}
    \Phi^{(n)}(h k, s) = \det(k_1)^{n_1} \det(k_2)^{n_2} \Phi^{(n)}(h)
    \end{equation}
for all $h \in H(\R)$ and $k \in K^{\circ}$, where
    \begin{equation}
    k = 
    \begin{pmatrix}
    a & b \\
    - b & a 
    \end{pmatrix} \in K_{v}
    \quad \quad
    k_1 = a + i b
    \quad \quad
    k_2 = a - i b.
    \end{equation}

If $y = a {}^t \overline{a}$ for some $a \in \GL_m(\C)$, a computation (omitted) shows
    \begin{equation}\label{equation:part_II:Eisenstein:setup:adelic_v_classical:explicit_Archimedean_scalar_weight_vector}
    \chi_v(a)^{-1} (\det y)^{- n/2} \Phi^{(n)}(w^{-1}n(b)m(a)) = (\det y)^{s - s_0} \det(-i y + b)^{-(s - s_0)} \det(iy + b)^{-(s - s_0) - n}
    \end{equation}
for any $b \in \mrm{Herm}_m(\R)$, where $s_0 = (n - m)/2$ (reduce to the case $a = 1_m$ and write $w^{-1} n(b) = n(-b(1_m + b^2)^{-1}) m(b + i 1_m)^{-1} k$ for $k \in K^{\circ}$).
Equation \eqref{equation:part_II:Eisenstein:setup:adelic_v_classical:explicit_Archimedean_scalar_weight_vector} may be used to translate various statements from \cite{Shimura82} to statements about Archimedean Whittaker functions, etc. (see Section \ref{ssec:Archimedean_identity:more_Whittaker} for more on this).

\begin{remark}\label{remark:part_II:Eisenstein:setup:adelic_v_classical:log_det_definition}
Given $g = x_g + i y_g \in M_{m,m}(\C)$ with $x_g, y_g$ Hermitian and $x_g$ positive definite, we define $\log \det(g)$ by the ``principal branch'' (such that $g \mapsto \log \det g$ is holomorphic, and $\log \det g \in \R$ if $y_g = 0$) as in \cite[{(1.11)}]{Shimura82} and the surrounding discussion of loc. cit.. If $y_g$ is positive definite and $x_g$ is only assumed Hermitian, we also take
    \begin{equation}
    \log \det g = \log \det (-i g) + m \log i \quad \quad \log \det \overline{g} = \log \det (i g) - m \log i
    \end{equation}
where $\log i \coloneqq \pi i/2$ (as in \cite[{(1.11)}]{Shimura82}).
This convention is implicit in \eqref{equation:part_II:Eisenstein:setup:adelic_v_classical:explicit_Archimedean_scalar_weight_vector}.
\end{remark}
    
            \subsection{Normalized local Whittaker functions}
            \label{ssec:part_II:local_Whittaker_Archimedean:normalized_Whittaker}
                Let $\chi \colon \C^{\times} \ra \C^{\times}$ and $\psi \colon \R \ra \C^{\times}$ be characters with $\psi(x) \coloneqq e^{2 \pi i x}$. Given a standard section $\Phi \in I(s, \chi)$, and given $T \in \mrm{Herm}_m(\R)$ with $\det T \neq 0$, there is a \emph{local Whittaker function} defined by the absolutely convergent integral
    \begin{equation}
    W_{T}(h,s,\Phi) \coloneqq \int_{N(\R)} \Phi(w^{-1} n(b) h, s) \psi(- \mrm{tr}(Tb)) ~dn(b)
    \end{equation}
for $h \in H(\R)$ and $s \in \C$ with $\Re(s) > m / 2$, where $dn(b)$ is the Haar measure with respect to the pairing $b,b' \mapsto \psi(\mrm{tr}(b b'))$ on $\mrm{Herm}_m(\R) \cong N(\R)$. For any fixed $h$, the function $W_T(h,s,\Phi)$ admits holomorphic admits holomorphic continuation to $s \in \C$, e.g. \cite[{\S 6}]{Ichino04}. Extending linearly defines $W_{T}(h, s, \Phi)$ whenever $\Phi$ is a finite meromorphic linear combination of standard sections. 

Next, fix $n \in \Z$ and assume $\chi|_{\R^{\times}} = \eta^n$. Set $s_0 \coloneqq (n - m)/2$.
We define a \emph{normalized local Whittaker function}
    \begin{align}
    W_{T}^*(h, s)^{\circ}_n \coloneqq \Lambda_{T}(s)^{\circ}_n W_{T}(h, s, \Phi^{(n)})
    \end{align}
using a \emph{local normalizing factor}, which we define as
    \begin{align}
    \Lambda_{T}(s)^{\circ}_n &\coloneqq \frac{(2 \pi)^{m(m-1)/2}}{(-2 \pi i)^{nm}} \pi^{m(- s + s_0)} \left ( \prod_{j = 0}^{m - 1} \Gamma(s - s_0 + n - j) \right ) |\det T|^{- s - s_0}_{F_0}.
    \label{equation:part_I:Eisenstein:local_Whittaker:Archimedean:normalizing_factor}
    \end{align}
This formula should be compared with \cite[{(3.3.14)}]{GS19} and \cite{Shimura82}.

We have a functional equation
    \begin{align}
    W_{T}^*(h, s)^{\circ}_n & = \eta(\det T)^{n - m - 1} W^*_T(h, s)^{\circ}_n. \label{equation:part_I:local_functional_equations:Archimedean:scalar_weight}
    \end{align}
The case when $T$ is positive definite follows from \cite[{Theorem 3.1}]{Shimura82} (via \crefext{IV:equation:Eisenstein:setup:adelic_v_classical:explicit_Archimedean_scalar_weight_vector}, see also \cite[{(3.54)}]{GS19}). The case of general $T$ (still with $\det T \neq 0$) should follow from \cite[{Theorem 4.2, (4.34.K)}]{Shimura82}, though we will give an alternative proof in our companion paper \crefext{IV:lemma:Eisenstein:local_functional_equations:Archimedean:scalar_weight}.
    
We use the notation
    \begin{align}\label{equation:part_I_only:local_Whittaker:normalized_Whittaker:only_s_var}
    W^*_T(s)^{\circ}_n \coloneqq W^*_T(1_{2m},s)^{\circ}_n e^{2 \pi \mrm{tr}(T)}
    \end{align}
and have the invariance property
    \begin{equation}\label{equation:part_I_only:local_Whittaker:normalized_Whittaker:Archimedean_k_invariance}
    W^*_T(s)^{\circ}_n = W^*_{{}^t \overline{k} T k}(s)^{\circ}_n
    \end{equation}
for any $k \in U(m)$ (unitary group for the positive definite Hermitian space specified by the Gram matrix $1_m$).

Next, write $(r_1,r_2)$ for the signature of $T$ (temporary notation). If either $n \geq r_1$ or $r_2 = 0$, then the function $W_{T}^*(h,s)^{\circ}_n$ is holomorphic for all $s \in \C$, for fixed $h$ (follows from \cite[{Theorem 4.2, (4.34.K)}]{Shimura82}).
We also have
    \begin{equation}\label{equation:part_I:Eisenstein:local_Whittaker:Archimedean:special_value}
    W^*_{T}(s_0)^{\circ}_n
    =
    \begin{cases}
    1 & \text{if $T$ is positive definite} \\
    0 & \text{if $m \leq n$ and $T$ is not positive definite.}
    \end{cases}
    \end{equation}
For the case when $T$ is positive definite, see \cite[{(3.15)}]{Shimura97} (also the proof of \cite[{Proposition 3.2}]{GS19}). The non positive definite case with $m \leq n$ follows from \cite[{Theorem 4.2, (4.34.K)}]{Shimura82} (see also \cite[{Proposition 3.3(i)}]{GS19}).


    \section{Archimedean local main theorem}
        \label{sec:Archimedean_identity}
                Let $V$ be the non-degenerate $\C / \R$ Hermitian space of rank $n$ and signature $(n - 1, 1)$.

We freely use notation for the Hermitian symmetric domain $\mc{D}$ and its special cycles (Section \ref{ssec:Hermitian_domain:setup}) as well as the Archimedean local Whittaker functions $W_{T}^*(s)^{\circ}_n$ for $T \in \mrm{Herm}_m(\R)$ (complex Hermitian matrices) with $\det T \neq 0$ \eqref{equation:part_I_only:local_Whittaker:normalized_Whittaker:only_s_var}.
            
            \subsection{Statement of local main theorem} 
            \label{ssec:Archimedean_identity:statement}
                Our main Archimedean local identity (``Archimedean local arithmetic Siegel--Weil'') is the following.

    \begin{theorem}\label{theorem:local_identities:Archimedean_identity:statement:main_Archimedean}
    Let $\underline{x} \in V^{m}$ be a $m$-tuple with nonsingular Gram matrix, and set $T = (\underline{x}, \underline{x})$. If $m \geq n - 1$ or if $T$ is not positive definite, we have
    \begin{equation}\label{equation:local_identities:Archimedean_identity:main_Archimedean}
    \int_{\mc{D}} [\xi(\underline{x})] \wedge c_1(\widehat{\mc{E}}^{\vee})^{n-m} = \frac{d}{d s} \bigg|_{s = -s_0} W^*_{T}(s)^{\circ}_n.
    \end{equation}
    where $s_0 = (n - m)/2$.
    \end{theorem}

In Theorem \ref{theorem:local_identities:Archimedean_identity:statement:main_Archimedean}, integration of the current $[\xi(\underline{x})] \wedge c_1(\widehat{\mc{E}}^{\vee})^{n - m}$ over $\mc{D}$ is understood in the sense discussed in \crefext{III:ssec:arith_cycle_classes:horizontal}. The displayed integral is convergent (combine Lemma \ref{lemma:convergence:integral_convergence} and Lemma \ref{lemma:convergence:KM_form_growth_bound}). The local functional equation \eqref{equation:part_I:local_functional_equations:Archimedean:scalar_weight} implies that the derivative of $W^*_{T}(s)^{\circ}_n$ at $s = s_0$ and $s = - s_0$ are the same up to a simple sign. 

The case $m = n$ of Theorem \ref{theorem:local_identities:Archimedean_identity:statement:main_Archimedean} is the content of \cite[{Theorem 4.17}]{Liu11} (when translating to Liu's notation, recall also that $W^*_{T}(0)^{\circ}_n = 0$ when $m \leq n$ for non positive-definite $T$ \eqref{equation:part_I:Eisenstein:local_Whittaker:Archimedean:special_value}. We do not give a new proof of this case. Indeed, we reduce the other cases of Theorem \ref{theorem:local_identities:Archimedean_identity:statement:main_Archimedean} to the case $m = n$ by the following limiting result.

\begin{proposition}\label{proposition:local_identities:Archimedean_identity:statement:limiting}
Let $T^{\flat} \in \mrm{Herm}_m(\R)$ be a matrix with $\det T^{\flat} \neq 0$, assume $m \leq n$, and set $s_0 = (n - m)/2$. Assume that either $m = n - 1$ or that $T^{\flat}$ is not positive definite. Given $t = \operatorname{diag}(t_{m+1}, \ldots, t_n) \in \mrm{Herm}_{n - m}(\R)$, set $T = \operatorname{diag}(t, T^{\flat})$.
    \begin{equation}
    \frac{d}{ds} \bigg |_{s = - s_0} W^*_{T^{\flat}}(s)^{\circ}_n = \lim_{t \ra 0^{\pm}} \left ( \frac{d}{ds} \bigg|_{s = 0} W^*_{T}(s)^{\circ}_n + (\log |t|_{F^+_v} + \log(4\pi) - \Gamma'(1)) W^*_{T^{\flat}}(-s_0)^{\circ}_n \right )
    \end{equation}
where $|t|_{\R} \coloneqq |\det t|_{\R}$, and where the sign on $0^{\pm}$ (meaning all $t_j$ have this sign) is
    \begin{equation}
    \begin{cases}
    - & \text{if $T^{\flat}$ is positive definite} \\
    + & \text{else.}
    \end{cases}
    \end{equation}
\end{proposition}

Recall that a similar limit formula was crucial for our proof of a ``non-Archimedean local Siegel--Weil formula'' in our companion paper \crefext{I:sec:non-Arch_identity}. A side-by-side comparison of the limit formulas in all cases (non-Archimedean and Archimedean) may be found in \crefext{I:ssec:part_I:local_Whittaker:limits}. Nevertheless, the proofs we give for these limit formulas are fairly different in the non-Archimedean vs Archimedean cases.

\begin{remark}\label{remark:local_identities:Archimedean_identity:statement:limiting}
In the situation of Proposition \ref{proposition:local_identities:Archimedean_identity:statement:limiting}, recall that $W^*_{T^{\flat}}(-s_0)^{\circ}_n = 0$ if $T^{\flat}$ is not positive definite, by \eqref{equation:part_I:Eisenstein:local_Whittaker:Archimedean:special_value} and the functional equation \eqref{equation:part_I:local_functional_equations:Archimedean:scalar_weight}.
Due to this vanishing, the term $(\log|t|_{F^+_v} + \log(4 \pi) - \Gamma'(1))$ from Proposition \ref{proposition:local_identities:Archimedean_identity:statement:limiting} should not be taken seriously outside the positive definite $T^{\flat}$ case (especially if $m \neq n - 1$).

If $T^{\flat}$ has signature $(p,q)$ for $q \geq 2$, we also have
    \begin{equation}\label{equation:Archimedean_identity:statement:higher_vanishing_order}
    \frac{d}{ds} \bigg|_{s = - s_0} W^*_{T^{\flat}}(s)^{\circ}_n = \frac{d}{ds} \bigg|_{s = 0} W^*_{T}(s)^{\circ}_n = 0
    \end{equation}
for any $t \in \mrm{Herm}_{n-m}(\R)$ with $\det t \neq 0$ by \cite[{Theorem 4.2, (4.34.K)}]{Shimura82}. Thus Proposition \ref{proposition:local_identities:Archimedean_identity:statement:limiting} holds for $T^{\flat}$ of signature $(p,q)$ when $q \geq 2$ (both sides of the identity are $0$).
\end{remark}

The proof of the remaining cases of Proposition \ref{proposition:local_identities:Archimedean_identity:statement:limiting} will occupy most of the rest of Section \ref{sec:Archimedean_identity}. The case of $T^{\flat}$ having signature $(m-1,1)$ is completed in Section \ref{ssec:Archimedean_identity:non_pos_def}, and the case of positive definite $T^{\flat}$ is completed in Section \ref{ssec:Archimedean_identity:pos_def}. We also obtain an explicit formula for both sides of \eqref{equation:local_identities:Archimedean_identity:main_Archimedean} when $T$ is positive definite, namely \eqref{equation:local_identities:Archimedean_identity:pos_def:explicit_formula_pos_def} (the formula is a polynomial in the eigenvalues of $T^{-1}$).

Once the proposition is proved, Theorem \ref{theorem:local_identities:Archimedean_identity:statement:main_Archimedean} follows (and is equivalent to the proposition for any given $T^{\flat}$, which is the $T$ in Theorem \ref{theorem:local_identities:Archimedean_identity:statement:main_Archimedean}) by the following argument.

\begin{proof}[Proof of equivalence of Theorem \ref{theorem:local_identities:Archimedean_identity:statement:main_Archimedean} and Proposition \ref{proposition:local_identities:Archimedean_identity:statement:limiting}]
Let $T^{\flat}$ be as in Proposition \ref{proposition:local_identities:Archimedean_identity:statement:limiting}. We may assume $T^{\flat}$ has signature $(m,0)$ or $(m-1,1)$ by Remark \ref{remark:local_identities:Archimedean_identity:statement:limiting}.
Suppose $\underline{x}^{\flat} = (x^{\flat}_1, \ldots, x^{\flat}_{m}) \in V^m$ satisfies $ (\underline{x}^{\flat}, \underline{x}^{\flat}) = T^{\flat}$. Given an orthogonal basis $\underline{x}^{\#} = (x_{m+1}, \ldots, x_{n})$ of $\operatorname{span}(\underline{x}^{\flat})^{\perp}$, set $t_j =  (x_j,x_j)$ for $j \geq m + 1$, set $t = (t_{m+1}, \ldots, t_n)$, set $\underline{x} = (x_{m+1}, \ldots, x_n, x^{\flat}_1, \ldots, x^{\flat}_{m}) \in V^n$, and set $T = (\underline{x}, \underline{x})$. We have
    \begin{equation}\label{equation:equivalence_limiting_decomposition}
    \frac{d}{d s} \bigg|_{s = 0} W^*_{T}(s)^{\circ}_n = \int_{\mc{D}} [\xi(\underline{x})] = \int_{\mc{D}} [\xi(\underline{x}^{\flat})] \wedge \omega(\underline{x}^{\#}) + \int_{\mc{D}(\underline{x}^{\flat})} \xi(x)
    \end{equation}
where the first equality is the $m = n$ case of Theorem \ref{theorem:local_identities:Archimedean_identity:statement:main_Archimedean} (already proved by Liu as cited above) and the second identity is by definition.

Using the pointwise convergence $\lim_{a \ra 0} \omega(a x) = c_1(\widehat{\mc{E}}^{\vee})$ on $\mc{D}$ for each $x \in V$ \eqref{equation:deform_KM_form}, we have
    \begin{equation}
    \lim_{\underline{x}^{\#} \ra 0} \int_{\mc{D}} [\xi(\underline{x}^{\flat})] \wedge \omega(\underline{x}^{\#}) = \int_{\mc{D}} [\xi(\underline{x}^{\flat})] \wedge c_1(\widehat{\mc{E}}^{\vee})^{n-m}
    \end{equation}
(say, where the limit runs over $\underline{x}^{\#} = (a_{m+1} x_{m+1}, \ldots, a_n x_n)$ as $a_j \ra 0$ for all $j$) by dominated convergence (applying estimate from the proof of Lemma \ref{lemma:convergence:KM_form_growth_bound} and convergence from Lemma \ref{lemma:convergence:integral_convergence}, particularly convergence of \eqref{equation:convergence:integral_convergence:dominating_integral}).

The closed submanifold $\mc{D}(\underline{x}^{\flat}) \subseteq \mc{D}$ is a single point if $T^{\flat}$ is positive definite (in which case we assumed $m = n - 1$), and is empty if $T^{\flat}$ is not positive definite. We thus have
    \begin{equation}
    \int_{\mc{D}(\underline{x}^{\flat})} \xi(x) = 
    \begin{cases}
    -\Ei(4 \pi (x_n,x_n)) & \text{if $T^{\flat}$ is positive definite} \\
    0 & \text{else}.
    \end{cases}
    \end{equation}

Recall asymptotics for the function $\Ei$ \eqref{equation:Ei_function_asymptotics} and recall $\Gamma'(1) = - \gamma$. Recall the special value formulas from \eqref{equation:part_I:Eisenstein:local_Whittaker:Archimedean:special_value} (and the functional equation \eqref{equation:part_I:local_functional_equations:Archimedean:scalar_weight}).
We substitute into \eqref{equation:equivalence_limiting_decomposition} to obtain
    \begin{equation}\label{equation:local_identities:Archimedean_identity:statement:limit_equivalence}
    \int_{\mc{D}} [\xi(\underline{x}^{\flat})] \wedge c_1(\widehat{\mc{E}}^{\vee})^{n - m} = \lim_{t \ra 0^{\pm}} \left ( \frac{d}{ds} \bigg|_{s = 0} W^*_{T}(s)^{\circ}_n + (\log |t|_{F^+_v} + \log(4\pi) - \Gamma'(1)) W^*_{T^{\flat}}(-s_0)^{\circ}_n \right )
    \end{equation}
(where the sign on $0^{\pm}$ is the sign of $t$, determined by the signature of $T^{\flat}$) which proves the claimed equivalence.
\end{proof}

\begin{remark}\label{remark:local_identities:Archimedean_identity:linear_invariance_Whittaker}
For any $T \in \mrm{Herm}_m(\R)$ with $\det T \neq 0$, recall that the (normalized) Archimedean local Whittaker function $W^*_{T}(s)^{\circ}_n$ satisfies a certain ``linear invariance'' property, i.e. the local Whittaker function is unchanged if we replace $T$ by ${}^t \overline{k} T k$ for any $k \in U(m)$ \eqref{equation:part_I_only:local_Whittaker:normalized_Whittaker:Archimedean_k_invariance}.
It is thus enough to prove Proposition \ref{proposition:local_identities:Archimedean_identity:statement:limiting} when $T^{\flat}$ is diagonal.
\end{remark}

\begin{remark}\label{remark:Archimedean_identity:statement:linear_invariance_geometric}
Using the linear invariance property for Whittaker functions, the limiting relation in \eqref{equation:local_identities:Archimedean_identity:statement:limit_equivalence} implies that the quantity
    \begin{equation}\label{equation:local_identities:Archimedean_identity:statement:linear_invariance}
    \int_{\mc{D}} [\xi(\underline{x}^{\flat})] \wedge c_1(\widehat{\mc{E}}^{\vee})^{n - m}
    \end{equation}
from Theorem \ref{theorem:local_identities:Archimedean_identity:statement:main_Archimedean}
is similarly linearly invariant (i.e. does not change if $\underline{x}^{\flat}$ is replaced by $\underline{x}^{\flat} \cdot k$ for any $k \in U(m)$, where $\underline{x}^{\flat}$ is viewed as a row vector of elements in $V$). Stated alternatively, we observe that the linear invariance result of Liu for $\int_{\mc{D}} [\xi(\underline{x}^{\flat})] \wedge c_1(\widehat{\mc{E}}^{\vee})^{n - m}$ when $m = n$ \cite[{Proposition 4.10}]{Liu11} can be used to prove the analogous linear invariance in our setting via limiting, even before we have proved Theorem \ref{theorem:local_identities:Archimedean_identity:statement:main_Archimedean} or Proposition \ref{proposition:local_identities:Archimedean_identity:statement:limiting}.
\end{remark}
        
            \subsection{Computation when \texorpdfstring{$n = 2$}{n = 2}} 
            \label{ssec:Archimedean_identity:case_n_is_2}
                Before proving Theorem \ref{theorem:local_identities:Archimedean_identity:statement:main_Archimedean} via Proposition \ref{proposition:local_identities:Archimedean_identity:statement:limiting} in the later sections, we check the $n = 2$ case of Theorem \ref{theorem:local_identities:Archimedean_identity:statement:main_Archimedean} by direct computation (the case $n = 1$ and $m \neq n$ is trivial as both sides of the identity are trivially $0$). The proof for general $n$ (which proceeds differently, not relying on the $n = 2$ computation) begins in Section \ref{ssec:Archimedean_identity:more_Whittaker} below.

Take $n = 2$ throughout Section \ref{ssec:Archimedean_identity:case_n_is_2}, and suppose $T \in \R$ is nonzero. By \cite[{(1.29) and (3.3)}]{Shimura82} (translation via \eqref{equation:part_II:Eisenstein:setup:adelic_v_classical:explicit_Archimedean_scalar_weight_vector}) and some rearranging, we have the formula
    \begin{align}
    W^*_{T}(s)^{\circ}_n & = \Gamma(s - 1/2)^{-1} |4 \pi T|^{s - 1/2} \int_0^{\infty} e^{-4 \pi T u} (u + 1)^{s + 1/2} u^{s - 3/2} ~du \label{equation:case_n_is_2:positive1}
    \\
    & = 1 + \Gamma(s - 1/2)^{-1} |4 \pi T|^{s - 1/2} \int_0^{\infty} e^{-4 \pi T u} ((u + 1)^{s + 1/2} - 1) u^{s - 3/2} ~du \label{equation:case_n_is_2:positive2}
    \end{align}
if $T > 0$, where the integrals in \eqref{equation:case_n_is_2:positive1} and \eqref{equation:case_n_is_2:positive2} are convergent for $\mrm{Re}(s) > 1/2$ and $\mrm{Re}(s) > - 1/2$ respectively. 
We similarly have
    \begin{align}
    W^*_{T}(s)^{\circ}_n & = \Gamma(s - 1/2)^{-1} |4 \pi T|^{s - 1/2} \int_1^{\infty} e^{4 \pi T u} (u - 1)^{s + 1/2} u^{s - 3/2} ~du \label{equation:case_n_is_2:negative}
    \end{align}
if $T < 0$, where displayed integral is convergent for $\mrm{Re}(s) > -3/2$.

\begin{proposition}\label{proposition:local_identities:Archimedean_identities:case_n_is_2}
Given any nonzero $T \in \R$ and any $x \in V$ with $T = (x, x)$, we have
    \begin{equation}
    \int_{\mc{D}} \xi(x) c_1(\widehat{\mc{E}}^{\vee}) = - \frac{d}{d s} \bigg|_{s = 1/2} W^*_{T}(s)^{\circ}_n = 
    \begin{cases}
    (-4 \pi T)^{-1} & \text{if $T > 0$} \\
    (4 \pi T)^{-1} e^{4 \pi T} - \Ei(4 \pi T) & \text{if $T < 0$.}
    \end{cases}
    \end{equation}
\end{proposition}

The preceding proposition (proved below) shows that Theorem \ref{theorem:local_identities:Archimedean_identity:statement:main_Archimedean} holds when $n = 2$ (the functional equation implies $-\frac{d}{ds}|_{s = 1/2} W^*_{T}(s)^{\circ}_n = \frac{d}{ds}|_{s = - 1/2} W^*_{T}(s)^{\circ}_n$).

\begin{lemma}\label{lemma:local_identities:Archimedean_identities:case_n_is_2:automorphic}
For any nonzero $T \in \R$, we have
    \begin{equation}
    -\frac{d}{d s} \bigg|_{s = 1/2} W^*_{T}(s)^{\circ}_n = 
    \begin{cases}
    (-4 \pi T)^{-1} & \text{if $T > 0$} \\
    (4 \pi T)^{-1} e^{4 \pi T} - \Ei(4 \pi T) & \text{if $T < 0$.}
    \end{cases}
    \end{equation}
\end{lemma}
\begin{proof}
Recall that $\Gamma(s)^{-1} = s + O(s^2)$ near $s = 0$.
The integrals in \eqref{equation:case_n_is_2:positive2} and \eqref{equation:case_n_is_2:negative} are convergent and holomorphic at $s = 1/2$. Directly evaluating the integrals at $s = 1/2$ gives the claimed formulas.
\end{proof}

\begin{lemma}\label{lemma:local_identities:Archimedean_identities:case_n_is_2:geometric}
With $x \in V$ and $T \in \R$ as in the statement of Proposition \ref{proposition:local_identities:Archimedean_identities:case_n_is_2}, we have
    \begin{equation}
    \int_{\mc{D}} \xi(x) c_1(\widehat{\mc{E}}^{\vee}) = 
    \begin{cases}
    (-4 \pi T)^{-1} & \text{if $T > 0$} \\
    (4 \pi T)^{-1} e^{4 \pi T} - \Ei(4 \pi T) & \text{if $T < 0$.}
    \end{cases}
    \end{equation}
\end{lemma}
\begin{proof}
By \eqref{equation:Chern_form_taut_bundle_2}, we have
    \begin{equation}
    c_1(\widehat{\mc{E}}^{\vee}) = \frac{1}{2 \pi i} \partial \overline{\partial} (\log R(x,z)) = \frac{1}{2 \pi i} \frac{d z \wedge d \overline{z}}{(1 - z \overline{z})^2}.
    \end{equation}

If $T > 0$, we have
    \begin{align}
    \int_{\mc{D}} \xi(x) c_1(\widehat{\mc{E}}^{\vee}) & = \int_{\mc{D}} \int_1^{\infty} e^{-4 \pi T u z \overline{z} (1 - z \overline{z})^{-1}} u^{-1} ~du ~\frac{1}{2 \pi i} \frac{dz \wedge d \overline{z}}{(1 - z \overline{z})^2}
    \\
    & = - 2 \int_{0}^1 \int_1^{\infty} e^{-4 \pi T u r^2 (1 - r^2)^{-1}} u^{-1} (1 - r^2)^{-2} r ~ du ~dr
    \\
    & = - \int_{1}^{\infty} \int_0^{\infty} e^{-4 \pi T v u} u^{-1} ~ dv ~du
    \\
    & = (- 4 \pi T)^{-1}
    \end{align}
via the change of variables $v = r^2 (1 - r^2)^{-1}$.

If $T < 0$, we have
    \begin{align}
    \int_{\mc{D}} \xi(x) c_1(\widehat{\mc{E}}^{\vee}) & = \int_{\mc{D}} \int_1^{\infty} e^{4 \pi T u (1 - z \overline{z})^{-1}} u^{-1} ~du \frac{1}{2 \pi i} \frac{dz \wedge d \overline{z}}{(1 - z \overline{z})^2}
    \\
    & = - 2 \int_{0}^1 \int_1^{\infty} e^{4 \pi T u (1 - r^2)^{-1}} u^{-1} (1 - r^2)^{-2} r ~ du ~dr
    \\
    & = - \int_{1}^{\infty} \int_{1}^{\infty} e^{4 \pi T v u} u^{-1} ~ dv ~ du
    \\
    & = (4 \pi T)^{-1} \int_{1}^{\infty} e^{4 \pi T u} u^{-2} ~du
    \end{align}
via the change of variables $v = (1 - r^2)^{-1}$. We also have
    \begin{equation}
    \int_{1}^{\infty} e^{4 \pi T u} u^{-2} ~du = e^{4 \pi T} - (4 \pi T) \Ei(4 \pi T)
    \end{equation}
via integration by parts.
\end{proof}

\begin{proof}[Proof of Proposition \ref{proposition:local_identities:Archimedean_identities:case_n_is_2}]
Already proved by direct computation in Lemmas \ref{lemma:local_identities:Archimedean_identities:case_n_is_2:automorphic} and \ref{lemma:local_identities:Archimedean_identities:case_n_is_2:geometric}.
\end{proof}

            \subsection{More on Archimedean local Whittaker functions} 
            \label{ssec:Archimedean_identity:more_Whittaker}
                We use some special functions studied by Shimura \cite{Shimura82} to describe the Archimedean Whittaker functions $W^*_{T}(s)^{\circ}_n$ from above. We allow arbitrary $n \in \Z$ for the moment.

We first recall Shimura's definitions. Given an integer $m \geq 0$, set 
    \begin{equation}
    \Gamma_m(s) = \pi^{m(m-1)/2} \prod_{k = 0}^{m - 1} \Gamma(s - k)
    \end{equation}
as in \cite[{(1.17.K)}]{Shimura82}, where $\Gamma$ is the usual gamma function. Given Hermitian matrices $h, h'$, the notation $h > h'$ will mean that $h - h'$ is positive definite. For
    \begin{align*}
    & \a, \b \in \C \quad \quad g \in \mrm{Herm}_m(\R)_{>0} \quad \quad h \in \mrm{Herm}_m(\R)
    \\
    & z \in \mc{H}' \coloneqq \{ z = x + i y \in M_{m,m}(\C) \text{ with } x, y \in \mrm{Herm}_m(\R) \text{ and } x >0 \}
    \end{align*}
we set
    \begin{align}
    \xi(g, h; \a, \b) & \coloneqq \int_{\mrm{Herm}_m(\R)} e^{-2 \pi i \operatorname{tr}(hx)} \det(x + i g)^{- \a} \det(x - i g)^{-\b} ~dx
    \\
    \eta(g, h; \a, \b) & \coloneqq \int_{\substack{\mrm{Herm}_m(\R) \\ x > h \\ x > - h}} e^{- \operatorname{tr}(gx)} \det(x + h)^{\a - m} \det(x - h)^{\b - m} ~dx
    \\
    \zeta_m(z; \a, \b) & \coloneqq \int_{\mrm{Herm}_m(\R)_{>0}} e^{-\tr(zx)} \det(x + 1_m)^{\a - m} \det(x)^{\b - m} ~ dx
    \\
    \omega_m(z; \a, \b) & \coloneqq \Gamma_m(\b)^{-1} \det(z)^{\b} \zeta_m(z; \a, \b) \label{equation:local_identities:Archimedean_identity:more_Whittaker:omega}
    \\
    \zeta_{p,q}(g; \a, \b) & \coloneqq e^{-\tr(g)/2} \int_{\substack{\mrm{Herm}_m(\R) \\ x + \diag(1_p,0) > 0 \\ x + \diag(0, 1_q) > 0}} e^{-\tr(gx)} \det(x + \diag(1_p,0))^{\a - m} \det(x + \diag(0, 1_q))^{\b - m}
    \end{align}
initially defined for $\operatorname{Re}(\a), \operatorname{Re}(\b) \gg 0$.  All implicit measures in the integrals are Euclidean. See Remark \ref{remark:part_II:Eisenstein:setup:adelic_v_classical:log_det_definition} for the log branch convention.

The special functions $\xi, \eta, \zeta_m, \omega_m, \zeta_{p,q}$ are copied from \cite[{(1.25), (1.26), (3.2), (3.6), (4.16)}]{Shimura82}, respectively. Formulas relating $\xi$ and $\eta$, relating $\eta$ and $\zeta_m$, and relating $\eta$ and $\zeta_{p,q}$ are given in \cite[{(1.29), (3.3), (4.18)}]{Shimura82}. These will be used implicitly in our computations below.

Recall that $\omega_m(z; \a, \b)$ admits holomorphic continuation to all $(z, \a, \b) \in \mc{H}' \times \C^2$ (by \cite[{Theorem 3.1}]{Shimura82}), and that $\Gamma_q(\a - p)^{-1} \Gamma_p(\b - q)^{-1} \zeta_{p,q}(g; \a, \b)$ admits holomorphic continuation to all $(\a, \b)$, for any $g$ (by \cite[{Theorem 4.2}]{Shimura82}). We also recall the special value formula
    \begin{equation}\label{equation:Archimedean_omega_constant_value}
    \omega_m(z; m, \b) = \omega_m(z; a, 0) = 1
    \end{equation}
for all $\a, \b \in \C$ \cite[{(3.15)}]{Shimura82}.

We will also use the differential operator $\Delta \coloneqq \det( \partial/ \partial z_{j,k})$ on the space of matrices $z = (z_{j,k})_{j,k} \in M_{m,m}(\C)$ as in \cite[{(3.10.II)}]{Shimura82} (also \cite[{(4-20)}]{Liu11}). For any $u \in \mrm{Herm}_m(\R)_{>0}$, with $u^{1/2}$ denoting its unique positive definite Hermitian square-root, we have the relation
    \begin{align}
    & (-1)^m \Delta( e^{-\operatorname{tr} u z} \det(u z)^{-\b} \omega_m(u^{1/2} z u^{1/2}; \a, \b) ) \notag
    \\
    & = e^{-\operatorname{tr} u z} \det(u z)^{-\b} \det(u) \omega_m(u^{1/2} z u^{1/2}; \a + 1, \b)
    \end{align}
where $\Delta$ is applied to the $z$ variable, and where both sides are evaluated at $z \in \mc{H}'$. The preceding formula is a slight variant of \cite[{(3.12)}]{Shimura82} and \cite[{(4-21)}]{Liu11} (and can be verified by similar reasoning). We will use this formula in its equivalent form
    \begin{align}
    & (-1)^m e^{\operatorname{tr} u z} \Delta( e^{-\operatorname{tr} u z} \det(z)^{-\b} \omega_m(u^{1/2} z u^{1/2}; \a, \b) ) \notag
    \\
    & = \Gamma_m(\b)^{-1} \det(u)^{\b + 1} \zeta_m(u^{1/2} z u^{1/2}; \a + 1, \b). \label{equation:differential_operator_alpha_plus_one}
    \end{align}

\begin{remark} \hfill
\begin{enumerate}[(1)]
\item The special function $\xi$ (which takes multiple arguments) should not be confused with the Green function from Section \ref{ssec:Hermitian_domain:local_cycles} (which takes one argument), as should be clear from context. The same applies to $\eta$ the special function (which takes multiple arguments) and $\eta$ the quadratic character (which takes one argument).
\item The definition of $\zeta_{p,q}$ in \cite[{(4.16)}]{Shimura82} has a running assumption that ``$g$ is diagonal'', but we can make the same definition without this diagonal assumption.
\item Liu also uses these functions but with slightly different normalizations \cite[{\S 4A}]{Liu11}. We follow Shimura's normalizations.
\end{enumerate}
\end{remark}  

Given $T \in \mrm{Herm}_m(\R)$ with $\det T \neq 0$, we set (non-standard)
    \begin{equation}
    W^*_T(\a, \b) \coloneqq e^{2 \pi \operatorname{tr} T} \frac{2^{m(m-1)} \pi^{-m \b}}{(-2 \pi i)^{m(\a - \b)}} \Gamma_m(\a) |\det T|^{-\a + m} \xi(1_m, T; \a, \b)
    \end{equation}
for $\a, \b \in \C$ initially defined for $\operatorname{Re}(\a), \operatorname{Re}(\b) \gg 0$. We have 
    \begin{equation}\label{equation:local_identities:Archimedean_identity:more_Whittaker:one_var_two_var_comparison}
    W^*_T(s)^{\circ}_n = W^*_T(\a, \b) \quad \quad \text{ when } \a = s - s_0 + n \text{ and } \b = s - s_0
    \end{equation}
where $s_0 = (n - m)/2$ (see \eqref{equation:part_II:Eisenstein:setup:adelic_v_classical:explicit_Archimedean_scalar_weight_vector}). We caution the reader against confusing $W^*_T(\a,\b)$ with the symbol $W^*_T(y,s)^{\circ}_n$ appearing in \cref{ssec:part_II:intro:Arch_local_ASW} (also \crefext{IV:part:part_IV:Eisenstein}).

For any $c \in \GL_m(\C)$ such that $c 1_{p,q} {}^t \overline{c} = T$ (where $1_{p,q} = \operatorname{diag}(1_p, -1_q) \in M_{m,m}(\R)$), and with $g \coloneqq {}^t \overline{c} c$, we have
    \begin{align}
     W^*_{T}(\a, \b) & = e^{2 \pi \operatorname{tr} T} (2 \pi)^{2 m \b} \pi^{-m \b} \Gamma_m(\b)^{-1} |\det T|^{-\a + m} 2^{m(m - \a - \b)} |\det T|^{\a + \b - m} \eta(2 \pi g, 1_{p,q}; \a, \b) \notag
     \\
     & = e^{2 \pi \operatorname{tr} T} \Gamma_m(\b)^{-1} |\det 4 \pi T|^{\b} \zeta_{p,q}(4 \pi g; \a, \b).
    \end{align}
When $T$ is positive definite, our conventions imply 
    \begin{equation}\label{equation:local_identities:Archimedean_identity:more_Whittaker:pos_def_W_T}
    W^*_T(\a, \b) = \omega_m(4 \pi g; \a, \b).    
    \end{equation}

\begin{lemma}\label{equation:local_identities:Archimedean_identity:more_Whittaker:alpha_deriv_vanishes}
Suppose $T \in \mrm{Herm}_m(\R)$ has $\det T \neq 0$. If $T$ is positive definite (resp. not positive definite), the function $W^*_{T}(\a, \b)$ admits holomorphic continuation to all $(\a,\b) \in \C^2$ (resp. for $\operatorname{Re}(\a) > m - 1$ and all $\b$). In this region, we have
    \begin{equation}
    \frac{\partial}{\partial \a} W^*_{T}(\a, \b) = 0
    \end{equation}
for $\b = 0$ (resp. $\b \in \Z_{\leq 0}$).
\end{lemma}
\begin{proof}
Let $T$ have signature $(p,q)$ and let $g$ be as above. The holomorphic continuation of $\Gamma_p(\b - q)^{-1} \Gamma_q(\a - p)^{-1} \zeta_{p,q}(4 \pi g; \a, \b)$ to all $(\a, \b) \in \C^2$ (as recalled above from \cite[{Theorem 4.2}]{Shimura82}) implies that $W^*_{T}(\a, \b)$ admits holomorphic continuation to the region claimed.

When $T$ is positive definite, \eqref{equation:Archimedean_omega_constant_value} implies $(\partial/\partial \a) W^*_T(\a, 0) = 0$. If $T$ is not positive definite, the function $\Gamma_m(\b)^{-1} \Gamma_p(\b - q)$ has a zero at every $\b \in \Z_{\leq 0}$, which implies $W^*_T(\a, \b) = 0$ for all $\b \in \Z_{\leq 0}$. Thus $(\partial/\partial \a) W^*_T(\a, \b) = 0$ for all $b \in \Z_{\leq 0}$ in this case.
\end{proof}

Suppose $n \geq 1$ is an integer. For any $g = (- 4 \pi)^{-1} \diag(a, b) \in \mrm{Herm}_n(\R)_{>0}$ with $a \in \mrm{Herm}_{n-1}(\R)_{<0}$ and $b \in \R_{<0}$, we have (as in \cite[{(4.25)}]{Shimura82} and also \cite[{(4-15)}]{Liu11})
    \begin{align}
    & e^{2 \pi \operatorname{tr} g} \zeta_{n-1, 1}(4 \pi g; \a, \b) 
    \\
    & = \int_{\C^{n-1}} e^{\operatorname{tr}(a w w^*) + b w^* w} \zeta_1(-b (1 + w^* w); \b, \a - n + 1) 
    \\
    & \mathrel{\hphantom{= \int_{\C^{n-1}}}} \cdot e^{\operatorname{tr}(-a u/2)} \eta(-a, u/2; \a, \b-1) ~dw \notag
    \\
    & = \int_{\C^{n-1}} e^{\operatorname{tr}(a w w^*) + b w^* w} \zeta_1(-b (1 + w^* w); \b, \a - n + 1) \label{equation:local_identities:Archimedean_identity:more_Whittaker:unfold_2}
    \\
    & \mathrel{\hphantom{= \int_{\C^{n-1}}}} \cdot \det(u)^{\a + \b - n} \zeta_{n-1}(- u^{1/2} a u^{1/2}; \a, \b-1) ~dw \notag  
    \end{align}
with $w \in \C^{n-1}$ viewed as column vectors, with $w^* \coloneqq {}^t \overline{w}$, with $u = 1_{n - 1} + w w^*$, with $u^{1/2}$ the unique positive definite Hermitian square-root of $u$, and with $dw$ being the Euclidean measure.

We next specialize \eqref{equation:local_identities:Archimedean_identity:more_Whittaker:unfold_2} to $\a = n$. We have
    \begin{equation}
    e^{b(1 + w^* w)} \zeta_1(-b(1 + w^* w); \b, 1) = \int_{1}^{\infty} e^{b(1 + w^* w) x} x^{\b - 1} ~dx.
    \end{equation}
Combining \eqref{equation:differential_operator_alpha_plus_one} and \eqref{equation:Archimedean_omega_constant_value}, we also find
    \begin{equation}\label{equation:local_identities:Archimedean_identity:more_Whittaker:zeta_via_Delta}
    \det(u)^{\b} \zeta_{n-1}(- u^{1/2} a u^{1/2}; n, \b-1) = (-1)^{n-1} e^{-\tr a u} \Delta|_{z = -a} (e^{-\tr uz} \det(z)^{-\b + 1}).
    \end{equation}
Hence, we have
    \begin{align}
    & (-1)^{n-1} e^b \Gamma_{m-1}(\b - 1)^{-1} e^{2 \pi \operatorname{tr} g} \zeta_{n-1, 1}(4 \pi g; n, \b) 
    \\
    & = \int_{\C^{n-1}} \int_1^{\infty} e^{\operatorname{tr}(a w w^*)} e^{b(1 + w^* w) x} x^{\b - 1} \label{equation:local_identities:Archimedean_identity:more_Whittaker:unfold_Delta}
    \\
    & \mathrel{\hphantom{= \int_{\C^{n-1}} \int_{1}^{\infty}}} \cdot e^{-\operatorname{tr} (1_{m-1} + w w^*) a} \Delta|_{z = -a} (e^{-\operatorname{tr} (1_{m-1} + w w^*) z} \det(z)^{-\b + 1}) ~dx ~dw. \notag
    \end{align}
These rearrangements are initially valid for $\operatorname{Re}(\b) \gg 0$, but in fact hold for all $\b \in \C$ by analytic continuation (see also \cite[{(3.8)}]{Shimura82} for estimates on $\zeta_1$ and $\zeta_{n-1}$ giving convergence).

The next lemma generalizes a calculation of Liu \cite[{Lemma 4.7}]{Liu11}, and will be used to re-express \eqref{equation:local_identities:Archimedean_identity:more_Whittaker:unfold_Delta} more explicitly. In the statement and proof below, we adopt the following notation from \cite[{Lemma 4.7}]{Liu11}: given a matrix $u \in M_{n,n}(\C)$ and sets $I,J \subseteq \{1, \ldots n\}$ of the same cardinality, the symbol $|u^{I,J}|$ (resp. $|u_{I,J}|$) will mean the determinant of the matrix obtained from $u$ by discarding (resp. keeping) the rows in indexed by $I$ and the columns indexed by $J$.

\begin{lemma}\label{lemma:Delta_operator_explicit_combinatorial}
Given any $u \in M_{m,m}(\C)$ and $z_0 \in \mrm{Herm}_m(\R)_{>0}$ with $z_0$ diagonal, we have
    \begin{equation}
    \Delta|_{z = z_0} (e^{\operatorname{tr}uz} \det(z)^s) = e^{\operatorname{tr} u z_0} \det(z_0)^{s} \sum_{t = 0}^m \sum_{\substack{J = \{j_1 < \cdots < j_t\} \\ J \subseteq \{1, \ldots, m\}}} \left( \prod_{k = 1}^t (s + k - 1) \right )|g_{0,J,J}|^{-1} |u^{J,J}|
    \end{equation}
for all $s \in \C$, where the inner sum runs over all subsets $J \subseteq \{1, \ldots, m\}$ of size $t$.
\end{lemma}
\begin{proof}
Observe that (upon fixing $u$ and $z_0$), the expression
    \begin{equation}
    e^{-\operatorname{tr} u z_0} \det(z_0)^{-s} \Delta|_{z = z_0} (e^{\operatorname{tr} u z} \det(z)^s)
    \end{equation}
is a polynomial in $s$. Hence it is enough to prove the lemma holds for all $s \in \Z_{\geq 1}$. The case $s = 1$ is given by the proof of \cite[{Lemma 4.7}]{Liu11} via combinatorial calculation. For all $s \in \Z_{\geq 1}$, a similar calculation shows
    \begin{equation}
    e^{- \operatorname{tr} u z_0} \det(z_0)^{-s} \Delta|_{z = z_0} (e^{\operatorname{tr}u z} \det(z)^s) = \sum_{t = 0}^m \sum_{\substack{J = \{j_1 < \cdots < j_t\} \\ J \subseteq \{1, \ldots, m\}}} N_{s,t} \cdot |z_{0,J,J}|^{-1} |u^{J,J}|
    \end{equation}
for all $s \in \Z_{\geq 1}$, where $N_{s,t}$ is the number of tuples $((\sigma_1, J_1), \ldots, (\sigma_s, J_s))$ where $J_i \subseteq J$ are disjoint subsets (possibly empty) with $\bigcup J_i = J$ and each $\sigma_i$ is a permutation of $J_i$. If $|J_i|$ denotes the cardinality of $J_i$, then there are $\binom{t + s - 1}{s - 1}$ possibilities for the tuple $(|J_1|, \ldots, |J_s|)$, and each such tuple admits $t!$ corresponding tuples $((\sigma_1, J_1), \ldots, (\sigma_s, J_s))$. 
Hence $N_{s,t} = \binom{t + s - 1}{s - 1} t! = \prod_{k = 1}^t (s + k - 1)$.
\end{proof}
        
            \subsection{Limiting identity: non positive definite \texorpdfstring{$T^{\flat}$}{Tflat}} 
            \label{ssec:Archimedean_identity:non_pos_def}
                Take integers $m, n \geq 1$, assume $m \leq n$, and set $s_0 = (n - m)/2$. Given $a = \diag(a_1, \ldots, a_{n-1}) \in \mrm{Herm}_{n-1}(\R)_{< 0}$ and $b \in \R_{<0}$, set $a^{\flat} = \diag(a_{n-m+1}, \ldots, a_{n-1}) \in \mrm{Herm}_{n-m}(\R)$ and
    \begin{align*}
    T &= (-4 \pi)^{-1} \diag(a, -b)  && T^{\flat} = (-4 \pi)^{-1} \diag (a^{\flat}, -b)
    \\
    g & = (-4 \pi)^{-1} \diag(a,b) && g^{\flat} = (-4 \pi)^{-1} \diag(a^{\flat}, b).
    \end{align*}
We have $T, g \in \mrm{Herm}_n(\R)$ and $T^{\flat}, g^{\flat} \in \mrm{Herm}_m(\R)$.

We have
    \begin{align*}
    W^*_{T^{\flat}}(\a, \b) &= e^{2 \pi \operatorname{tr} T^{\flat}} \Gamma_m(\b)^{-1} |\det 4 \pi T^{\flat}|^{\b} \zeta_{m-1,1}(4 \pi g^{\flat}; \a, \b)
    \\
    & = e^{b} |\det 4 \pi T^{\flat}|^{\b} \pi^{- m + 1} \Gamma(\b)^{-1} \Gamma_{m-1}(\b-1)^{-1} e^{2 \pi \operatorname{tr} g^{\flat}} \zeta_{m - 1,1}(4 \pi g^{\flat}; \a, \b)
    \end{align*}
which implies
    \begin{align}
    & \frac{\partial}{\partial \b} W^*_{T^{\flat}}(m, \b) \label{equation:local_identities:Archimedean_identity:non_pos_def:partial_beta}
    \\
    & = \left ( \frac{d}{d \beta} \Gamma(\b)^{-1} \right ) |\det 4 \pi T^{\flat}|^{\b} \pi^{- m + 1} e^{b} \Gamma_{m - 1}(\b - 1)^{-1} e^{2 \pi \operatorname{tr} g^{\flat}} \zeta_{m - 1,1}(4 \pi g^{\flat}; m, \b) \notag
    \end{align}
whenever both sides are evaluated at $\b \in \Z_{\leq 0}$.

Equation \eqref{equation:local_identities:Archimedean_identity:more_Whittaker:one_var_two_var_comparison} and Lemma \ref{equation:local_identities:Archimedean_identity:more_Whittaker:alpha_deriv_vanishes} imply
    \begin{equation}
    \frac{d}{ds} \bigg|_{s = - s_0} W^*_{T^{\flat}} (s)^{\circ}_n = \frac{\partial}{\partial \b} \bigg|_{\b = m - n} W^*_{T^{\flat}}(m, \b).
    \end{equation}
Since $\Gamma(s)^{-1}$ has residue $(-1)^{n - m} (n - m)!$ at $s = m - n$, we use \eqref{equation:local_identities:Archimedean_identity:non_pos_def:partial_beta} and \eqref{equation:local_identities:Archimedean_identity:more_Whittaker:unfold_Delta} to find
    \begin{align}
    & \frac{d}{ds} \bigg|_{s = - s_0} W^*_{T^{\flat}} (s)^{\circ}_n
    \\
    & = (-1)^{n - m} (n - m)! |\det 4 \pi T^{\flat}|^{m - n} (-\pi)^{-m + 1}
    \\
    & \mathrel{\hphantom{=}} \cdot \int_{\C^{m-1}} \int_1^{\infty} e^{\operatorname{tr}(a w w^*)} e^{b(1 + w^* w) x} x^{m - n - 1} \notag
    \\
    & \mathrel{\hphantom{= \cdot \int_{\C^{n-1}} \int_{1}^{\infty}}} \cdot e^{-\operatorname{tr} (1_{m-1} + w w^*) a^{\flat}} \Delta|_{z = -a^{\flat}} (e^{-\operatorname{tr} (1_{m-1} + w w^*) z} \det(z)^{n - m + 1}) ~dx ~dw. \notag
    \end{align}
Next, we write $w = (w_1, \ldots, w_m)$ and apply Lemma \ref{lemma:Delta_operator_explicit_combinatorial} to find (using $\det(1 + w w^*) = 1 + w^* w$ as in \cite[{Lemma 2.2}]{Shimura82})
    \begin{align}
    & \frac{d}{ds} \bigg|_{s = - s_0} W^*_{T^{\flat}} (s)^{\circ}_n \notag 
    \\
    & = (-1)^{n + 1} (n - m)! |\det 4 \pi T^{\flat}|^{m - n} \pi^{-m + 1} \det(- a^{\flat})^{n - m} \notag
    \\
    & \mathrel{\hphantom{=}} \cdot \sum_{t = 0}^{m - 1} \sum_{\substack{I = \{i_1 < \cdots < i_t\} \\ I \subseteq \{1, \ldots, m - 1\}}} \bigg ( \prod_{k = 1}^{m - 1 - t} (n - m + k) (a_{i_1} \cdots a_{i_t}) \notag
    \\
    & \mathrel{\hphantom{= \cdot}} \int_{\C^{m-1}} \int_1^{\infty} e^{\operatorname{tr}(a w w^*)} e^{b(1 + w^* w) x} x^{m - n - 1} (1 + w_{i_1} \overline{w}_{i_1} + \cdots + w_{i_t} + \overline{w}_{i_t}) ~dx ~dw \bigg ) \notag
    \\
    & = (-1)^{n + 1} \pi^{-m + 1} (-b)^{m - n} \label{equation:local_identities:Archimedean_identity:non_pos_def:unfolded}
    \\
    & \mathrel{\hphantom{=}} \cdot \sum_{t = 0}^{m - 1} \sum_{\substack{I = \{i_1 < \cdots < i_t\} \\ I \subseteq \{1, \ldots, m - 1\}}} \bigg ( (n - 1 - t)! (a_{i_1} \cdots a_{i_t}) \notag
    \\
    & \mathrel{\hphantom{= \cdot}} \int_{\C^{m-1}} \int_1^{\infty} e^{\operatorname{tr}(a w w^*)} e^{b(1 + w^* w) x} x^{m - n - 1} (1 + w_{i_1} \overline{w}_{i_1} + \cdots + w_{i_t} + \overline{w}_{i_t}) ~dx ~dw \bigg ). \notag
    \end{align}
We have used exponential decay of the function $\int_1^{\infty} e^{c x} x^{m - n - 1}$ as $c \ra - \infty$ for convergence estimates (to rearrange integrals). The previous formulas also hold when $T^{\flat}, m, g^{\flat}, a^{\flat}$ are replaced by $T, n, g, a$ (the latter is just the special case $m = n$).

For the reader's convenience, we recall the formulas (which will be used below)
    \begin{equation}\label{equation:local_identities:Archimedean_identity:non_pos_def:reader_convenience_integrals}
    \int_{\R^2} e^{c(x^2 + y^2)} ~dx ~dy = - \pi c^{-1} \quad \quad \int_{\R^2} (x^2 + y^2) e^{c(x^2 + y^2)} ~dx ~dy = \pi c^{-2}
    \end{equation}
valid for any $c \in \R_{<0}$.

\begin{proof}[Proof of Proposition \ref{proposition:local_identities:Archimedean_identity:statement:limiting} when $T^{\flat}$ is not positive definite]
It is enough to check the case where $T^{\flat}$ is diagonal and signature $(m-1,1)$, by Remarks \ref{remark:local_identities:Archimedean_identity:statement:limiting} and \ref{remark:local_identities:Archimedean_identity:linear_invariance_Whittaker}. Take notation as above. There is nothing to check when $m = n$. Otherwise, we may show
    \begin{equation}
    \lim_{\substack{a_i \ra 0 \\ i = 1, \ldots, {n - m}}} \frac{d}{ds} \bigg|_{s = 0} W^*_{T}(s)^{\circ}_n = \frac{d}{ds} \bigg|_{s = - s_0} W^*_{T^{\flat}}(s)^{\circ}_n
    \end{equation}
via \eqref{equation:local_identities:Archimedean_identity:non_pos_def:unfolded}. Indeed, interchanging the limit and integrals (dominated convergence) and integrating out the variables $w_1, \ldots, w_{n - m}$ gives the claim (using the left identity in \eqref{equation:local_identities:Archimedean_identity:non_pos_def:reader_convenience_integrals}).
\end{proof}
        
            \subsection{Limiting identity: positive definite \texorpdfstring{$T^{\flat}$}{Tflat}} 
            \label{ssec:Archimedean_identity:pos_def}
                Take any integer $n \geq 1$ and set $m = n - 1$, so that $s_0 = (n - m)/2 = 1/2$. 
Given $a = \diag(a_1, \ldots, a_{n-1}) \in \mrm{Herm}_{n-1}(\R)_{< 0}$ and $b \in \R_{<0}$, set
    \begin{equation}
    T^{\flat} = (-4 \pi)^{-1} a \quad \text{and} \quad T = (- 4 \pi)^{-1} \diag(a, -b).
    \end{equation}

Equation \eqref{equation:local_identities:Archimedean_identity:more_Whittaker:one_var_two_var_comparison} and Lemma \ref{equation:local_identities:Archimedean_identity:more_Whittaker:alpha_deriv_vanishes} imply
    \begin{equation}\label{equation:local_identities:Archimedean_identity:pos_def:one_var_two_var}
    \frac{d}{ds} \bigg|_{s = - 1/2} W^*_{T^{\flat}} (s)^{\circ}_n = - \frac{d}{ds} \bigg|_{s = 1/2} W^*_{T^{\flat}} (s)^{\circ}_n = - \frac{\partial}{\partial \b} \bigg|_{\b = 0} W^*_{T^{\flat}}(n, \b)
    \end{equation}
where the first equality is via the functional equation \eqref{equation:part_I:local_functional_equations:Archimedean:scalar_weight}.
We have
    \begin{equation}
    W^*_{T^{\flat}}(n, \b) = \Gamma_m(\b)^{-1} \det(-a)^{\b} \zeta_m(-a; n, \b) = (-1)^m e^{- \tr a} \det(-a)^{-1} \Delta|_{z = 1_m} (e^{\tr a z} \det(z)^{-\b})
    \end{equation}
where the first equality is by \eqref{equation:local_identities:Archimedean_identity:more_Whittaker:omega} and \eqref{equation:local_identities:Archimedean_identity:more_Whittaker:pos_def_W_T}, and the second equality is by 
\eqref{equation:differential_operator_alpha_plus_one} and \eqref{equation:Archimedean_omega_constant_value}. Applying Lemma \ref{lemma:Delta_operator_explicit_combinatorial} then yields
    \begin{equation}\label{equation:local_identities:Archimedean_identity:pos_def:explicit_formula_pos_def}
    - \frac{\partial}{\partial \b} \bigg|_{\b = 0} W^*_{T^{\flat}}(n, \b) = \det(a)^{-1} \sum_{t = 0}^{m - 1} \sum_{\substack{I = \{i_1 < \cdots < i_t\} \\ I \subseteq \{1, \ldots, m\}}} (m - 1 - t)! (a_{i_1} \cdots a_{i_t}).
    \end{equation}

Before proceeding, we define several functions which serve only to aid computation in Section \ref{ssec:Archimedean_identity:pos_def}. Set
    \begin{align}
    d_m & \coloneqq \sum_{t = 0}^{m} \sum_{\substack{I = \{i_1 < \cdots < i_t\} \\ I \subseteq \{1, \ldots, m\}}} (m - t)! (a_{i_1} \cdots a_{i_t})
    \\
    q_m(x) & \coloneqq (x + a_1)^{-1} \cdots (x + a_m)^{-1}
    \\
    r_m(x) & \coloneqq 1 - (x + a_1)^{-1} - \cdots - (x + a_m)^{-1}
    \\
    h_m(x) & \coloneqq q(x) e^x \sum_{k = 0}^{m - 1} \sum_{t = 0}^k \sum_{\substack{I = \{i_1 < \cdots < i_t\} \\ I \subseteq \{1, \ldots, m\}}} (m - 1 - k)! a_{i_1} \cdots a_{i_t} x^{k - t}
    \\
    u_m(x) & \coloneqq \sum_{t = 0}^m \sum_{\substack{I = \{i_1 < \cdots < i_t\} \\ I \subseteq \{1, \ldots, m\}}} (m - t)! a_{i_1} \cdots a_{i_t} (1 - (x + a_{i_1})^{-1} - \cdots - (x + a_{i_t})^{-1})
    \\
    f_m(x) & \coloneqq q_m(x) u_m(x)
    \end{align}
where dependence on $a_1, \ldots, a_m$ is suppressed from notation.

Next, we consider \eqref{equation:local_identities:Archimedean_identity:non_pos_def:unfolded} for the matrix $T$. Changing variables $x \mapsto x/b$ and computing the $dw$ integral (using \eqref{equation:local_identities:Archimedean_identity:non_pos_def:reader_convenience_integrals}), we find
    \begin{equation}\label{equation:local_identities:Archimedean_identity:pos_def:separate_out_Ei}
    \frac{d}{ds} \bigg|_{s = 0} W^*_{T}(s)^{\circ}_n = - \int_{-\infty}^{1/b} f_m(x) e^x x^{-1} ~dx = -\Ei(b) + \int_{-\infty}^{1/b} (1 - f_m(x)) e^x x^{-1} ~dx
    \end{equation}
with $m = n - 1$ as above, and where $\Ei$ is the exponential integral function from Section \ref{ssec:Hermitian_domain:local_cycles}.

\begin{lemma}\label{lemma:local_identities:Archimedean_identity:pos_def:singularity}
We have $f_m(x) = 1 + O(x)$ near $x = 0$.
\end{lemma}
\begin{proof}
In the lemma statement, the variables $a_1, \ldots, a_m$ are understood to be fixed (and negative). Since $f_m(x)$ is a rational function of $x$, it is enough to check $f_m(0) = 1$, i.e. that
    \begin{align*}
    \sum_{t = 0}^m \sum_{\substack{I = \{i_1 < \cdots < i_t\} \\ I \subseteq \{1, \ldots, m\}}} (m - t)! a_{j_1}^{-1} \cdots a_{j_{m - t}}^{-1} (1 - a_{i_1}^{-1} - \cdots - a_{i_t}^{-1}) = 1
    \end{align*}
where $\{j_1, \ldots, j_{m - t} \} = \{1, \ldots, m\} \setminus \{i_1, \ldots, i_t\}$. This holds because the sum telescopes, i.e. for any given $t = 0, \ldots, m - 1$, we have
    \begin{equation*}
    \sum_{\substack{I = \{i_1 < \cdots < i_t\} \\ I \subseteq \{1, \ldots, m\}}} (m - t)! a_{j_1}^{-1} \cdots a_{j_{m - t}}^{-1} = \sum_{\substack{I' = \{i'_1 < \cdots < i'_{t+1}\} \\ I' \subseteq \{1, \ldots, m\}}} (m - t - 1)! a_{j'_1}^{-1} \cdots a_{j'_{m - t - 1}}^{-1} (a_{i'_1}^{-1} + \cdots + a_{i'_t}^{-1})
    \end{equation*}
where $\{j'_1, \ldots, j'_{m - t - 1}\} = \{1, \ldots, m\} \setminus \{i'_1, \ldots, i'_{t+1}\}$.
\end{proof}

\begin{lemma}\label{lemma:local_identities:Archimedean_identity:pos_def:antiderivative}
We have $\frac{d}{dx} h_m(x) = (1 - f_m(x)) e^x x^{-1}$.
\end{lemma}
\begin{proof}
We prove this by induction on $m$. The case $m = 0$ is clear, as both sides of the identity are $0$.
Next, suppose the claim holds for some $m$. We write $f_{m+1}(x)$, $h_{m+1}(x)$, etc. for the corresponding functions formed with respect to the tuple $(a_1, \ldots, a_m, a_{m+1})$ for any given choice of $a_{m+1} \in \R_{<0}$. Observe that we have an inductive formula
    \begin{equation}
    h_{m+1}(x) = h_m(x) + q_{m+1}(x) d_m
    \end{equation}
which implies
    \begin{equation}
    \frac{d}{dx} h_{m+1}(x) - \frac{d}{dx} h_m(x) = q_{m+1}(x) r_{m+1}(x) e^x d_m. 
    \end{equation}
So it is enough to check $f_m(x) - f_{m+1}(x) = x q_{m+1}(x) r_{m+1}(x) d_m$, which is equivalent to checking
    \begin{equation}\label{equation:local_identities:Archimedean_identity:pos_def:antiderivative_3}
    (x + a_{m+1}) u_m(x) - u_{m+1}(x) = x r_{m+1}(x) d_m.
    \end{equation}
To see that this holds, we first compute
    \begin{align*}
    & u_{m+1} - (a_{m+1} u_m(x) - a_{m+1}(x + a_{m+1})^{-1} d_m)
    \\
    & = \sum_{t = 0}^m \sum_{\substack{I = \{i_1 < \cdots < i_t\} \\ I \subseteq \{1, \ldots, m\}}} (m + 1 - t)! a_{i_1} \cdots a_{i_t} (1 - (x + a_{i_1})^{-1} - \cdots - (x + a_{i_t})^{-1}).
    \end{align*}
Using the identity $a_{m+1}(x + a_{m+1})^{-1} = 1 - x (x + a_{m+1})^{-1}$, we see that \eqref{equation:local_identities:Archimedean_identity:pos_def:antiderivative_3} is equivalent to the identity
    \begin{align}
    & x u_m(x) - x r_{m+1}(x) d_m + (1 - x (x + a_{m+1})^{-1}) d_m \label{equation:local_identities:Archimedean_identity:pos_def:antiderivative_6}
    \\
    & = \sum_{t = 0}^m \sum_{\substack{ I = \{i_1 < \cdots < i_t\} \\ I \subseteq \{1, \ldots, m\}}} (m + 1 - t)! a_{i_1} \cdots a_{i_t} (1 - (x + a_{i_1})^{-1} - \cdots - (x + a_{i_t})^{-1}). \notag
    \end{align}
To see that the latter identity holds, we compute
    \begin{align*}
    & x u_m(x) - x r_{m+1}(x) d_m + (1 - x (x + a_{m+1})^{-1}) d_m = x u_m(x) - x r_m(x) + d_m
    \\
    & = \sum_{t = 0}^m \sum_{\substack{I = \{i_1 < \cdots < i_t\} \\ I \subseteq \{1, \ldots, m\}}} (m - t)! a_{i_1} \cdots a_{i_t} (m + 1 - t - a_{j_1}(x + a_{j_1})^{-1} - \cdots - a_{j_{m-t}} (x + a_{j_{m-t}})^{-1})
    \\
    & = \sum_{t = 0}^m \sum_{\substack{I = \{i_1 < \cdots < i_t\} \\ I \subseteq \{1, \ldots, m\}}} (m + 1 - t)! a_{i_1} \cdots a_{i_t}
    \\
    & \mathrel{\hphantom{=}} - \sum_{t = 0}^m \sum_{\substack{I = \{i_1 < \cdots < i_t\} \\ I \subseteq \{1, \ldots, m\}}} (m - t)! a_{i_1} \cdots a_{i_t} (a_{j_1}(x + a_{j_1})^{-1} + \cdots + a_{j_{m-t}} (x + a_{j_{m-t}})^{-1})
    \end{align*}
where $\{j_1, \ldots, j_{m - t}\} = \{1, \ldots, m\} \setminus \{i_1, \ldots, i_t\}$.

We thus find that \eqref{equation:local_identities:Archimedean_identity:pos_def:antiderivative_6} is equivalent to the identity
    \begin{align*}
    & \sum_{t = 0}^m \sum_{\substack{I = \{i_1 < \cdots < i_t\} \\ I \subseteq \{1, \ldots, m\}}} (m + 1 - t)! a_{i_1} \cdots a_{i_t} ((x + a_{i_1})^{-1} + \cdots + (x + a_{i_t})^{-1})
    \\
    & = \sum_{t = 0}^m \sum_{\substack{I = \{i_1 < \cdots < i_t\} \\ I \subseteq \{1, \ldots, m\}}} (m - t)! a_{i_1} \cdots a_{i_t} (a_{j_1}(x + a_{j_1})^{-1} + \cdots + a_{j_{m-t}} (x + a_{j_{m-t}})^{-1})
    \end{align*}
with $\{j_1, \ldots, j_{m - t}\}$ as above, and this identity holds because both expressions are equal to
    \begin{equation*}
    \sum_{t = 0}^m \sum_{\substack{I = \{i_1 < \cdots < i_t\} \\ I \subseteq \{1, \ldots, m\}}} \sum_{i = 1}^t (m + 1 - t)! a_{i_1} \cdots a_{i_t} (x + a_{i_i})^{-1}. \qedhere
    \end{equation*}
\end{proof}

\begin{proof}[Proof of Proposition \ref{proposition:local_identities:Archimedean_identity:statement:limiting} when $T^{\flat}$ is positive definite]
We may assume $T^{\flat}$ is diagonal by Remark \ref{remark:local_identities:Archimedean_identity:linear_invariance_Whittaker}. With $T$ and $T^{\flat}$ as above, we find
    \begin{equation}
    \lim_{b \ra 0^-} \left ( \frac{d}{ds} \bigg|_{s = 0} W^*_{T}(s)^{\circ}_n + \Ei(b) \right ) = \int_{-\infty}^{0} (1 - f_m(x)) e^x x^{-1} ~dx
    \end{equation}
via \eqref{equation:local_identities:Archimedean_identity:pos_def:separate_out_Ei} (and Lemma \ref{lemma:local_identities:Archimedean_identity:pos_def:singularity} for convergence of the integral). The asymptotics for $\Ei(b)$ as $b \ra 0^-$ \eqref{equation:Ei_function_asymptotics} show that it is enough to verify the identity
    \begin{equation}
    \frac{d}{ds}\bigg|_{s = -1/2} W^*_{T^{\flat}}(s)^{\circ}_n = \int_{-\infty}^{0} (1 - f_m(x)) e^x x^{-1} ~dx.
    \end{equation}
The left-hand side was computed in \eqref{equation:local_identities:Archimedean_identity:pos_def:explicit_formula_pos_def} (via \eqref{equation:local_identities:Archimedean_identity:pos_def:one_var_two_var}). The right-hand side is equal to $h_m(0)$ (in the notation above) via the explicit antiderivative result from Lemma \ref{lemma:local_identities:Archimedean_identity:pos_def:antiderivative}. Inspecting the formula for $h_m(x)$ shows that the claimed identity holds.
\end{proof}       

    \clearpage
    

    \phantomsection
    \addcontentsline{toc}{part}{References}
    \renewcommand{\addcontentsline}[3]{}
    \printbibliography

\end{document}